\documentclass[%aip,
jmp,%
rsi,%
amsmath,amssymb,
preprint,%reprint,%
longbibliography
]{revtex4-1}

\usepackage{graphicx}% Include figure files
\usepackage{dcolumn}% Align table columns on decimal point
\usepackage{bm}% bold math
\usepackage[mathlines]{lineno}% Enable numbering of text and display math
\usepackage{wrapfig}
\usepackage{epstopdf}
\usepackage{epsfig}
\usepackage{calrsfs}
\usepackage{IEEEtrantools}
\usepackage{amsmath}
\usepackage{dcolumn}% Align table columns on decimal point
\usepackage{caption}
\usepackage{subcaption}
\usepackage{commath}
\usepackage{float}
\usepackage{rotating}
\usepackage{ textcomp }

\DeclareMathOperator{\sech}{sech}
\usepackage{mathtools}
\usepackage{dsfont}
\DeclarePairedDelimiter{\ceil}{\lceil}{\rceil}

\newcommand{\p}{\partial}
\newcommand{\Sch}{Schr\"{o}dinger}
\newcommand{\GL}{Gr\"{u}nwald–-Letnikov}
\newcommand{\RL}{Riemann--Liouville}
\newcommand{\Cp}{Caputo}
\newcommand{\ML}{Mittag-Leffler}

\newcommand{\nn}{\nonumber}
\newcommand{\cpg}{\clearpage}

\begin{document}

\preprint{AIP}

\title{Expansion of fractional derivatives in terms of an integer derivative series:
physical and numerical applications}

\author{Anastasia Gladkina$^{1}$}
 \email{agladkin@mines.edu}
\author{Gavriil Shchedrin$^{1}$}
% \email{shchedrin@mines.edu}
\author{U. Al Khawaja$^{2}$}
% \email{u.alkhawaja@uaeu.ac.ae}
\author{Lincoln D. Carr$^{1}$}
% \email{lcarr@mines.edu}
\affiliation{%
$^{1}$Physics Department, Colorado School of Mines, Golden, Colorado 80401, USA
}%
\affiliation{
$^{2}$Physics Department, United Arab Emirates University, P.O. Box 15551, Al-Ain, United Arab Emirates}

\date{\today}
\begin{abstract}
We use the displacement operator to derive an infinite series of integer order derivatives for the Gr\"{u}nwald-Letnikov fractional derivative and show its correspondence to the Riemann-Liouville and Caputo fractional derivatives.
We demonstrate that all three definitions of a fractional derivative lead to the same infinite series of integer order derivatives.
We find that functions normally represented by Taylor series with a finite radius of convergence have a corresponding integer derivative expansion with an infinite radius of convergence.
Specifically, we demonstrate robust convergence of the integer derivative series for the hyperbolic secant (tangent) function, characterized by a finite radius of convergence of the Taylor series $R=\pi/2$, which describes bright (dark) soliton propagation in non-linear media.
We also show that for a plane wave, which has a Taylor series with an infinite radius of convergence, as the number of terms in the integer derivative expansion increases, the truncation error decreases.
Finally, we illustrate the utility of the truncated integer derivative series by solving two linear fractional differential equations, where the fractional derivative is replaced by an integer derivative series up to the second order derivative.
We find that our numerical results closely approximate the exact solutions given by the Mittag-Leffler and Fox-Wright functions.
Thus, we demonstrate that the truncated expansion is a powerful method for solving linear fractional differential equations, such as the fractional {\Sch} equation.

%
%Valid PACS numbers may be entered using the \verb+\pacs{#1}+ command.
\end{abstract}

\pacs{Valid PACS appear here}% PACS, the Physics and Astronomy
                             % Classification Scheme.
%Use showkeys class option if keyword
                              %display desired
\maketitle

\section{Introduction}

%Task: please add more examples

%Task: add references after each example that you do not prove or derive in the paper.

Fractional calculus is a powerful tool to describe physical systems characterized by multiple time and length scales, nonlocality, fractional geometry, non-Gaussian statistics, 
and non-Fickian transport \cite{kilbas2006theory, samko1993fractional}. Anomalous diffusion through disordered media \cite{bouchaud1990anomalous, havlin1987diffusion}, 
hydrogeologic treatment of water propagation through soil and rocks \cite{
benson2013fractional, schumer2003multiscaling, schumer2001eulerian}, L\'{e}vy flights \cite{del2003front}, and turbulence \cite{kim1985application} are among the physical phenomena that can be consistently described within the framework of fractional calculus \cite{hilfer2000applications, west2014colloquium, herrmann2014fractional}. Similarly, certain biological systems, e.g., neuron clusters and heart cell arrays, exhibit multiple time scales that define fractional dynamics of the biological response to external stimuli \cite{lundstrom2008fractional, anastasio1994fractional}.

The building block of fractional calculus is a fractional derivative. There are multiple ways to generalize an integer order derivative to fractional order, and in this paper we exclusively concentrate on the {\RL}, {\Cp}, and {\GL} definitions \cite{kilbas2006theory}. The Riemann-Liouville and Caputo definitions are integral forms of the fractional derivative, especially suitable for solving linear fractional differential equations (FDEs) \cite{kilbas2006theory, samko1993fractional}. The Gr\"{u}nwald-Letnikov derivative is a discrete form of the fractional derivative, represented by a function summed over its history, and is primarily used in numerical methods to solve linear FDEs. The {\GL} derivative gives a computationally straightforward way to find the fractional derivative of an arbitrary function, yet it provides no direction towards finding its explicit analytical form. Except for a few trivial cases, where a fractional derivative can be expressed in terms of elementary or special functions, the {\RL} and {\Cp} derivatives also lead to expressions that are implicit or indirect
\cite{kilbas2006theory, samko1993fractional, herrmann2014fractional}.

%Moreover, we derive a formal bound on the truncation error in the integer order derivative expansion, thus, establishing a firm ground for its robust numerical applications of solving fractional FDEs.

Despite the fact that {\RL}, {\Cp}, and {\GL} definitions are three different forms of the fractional derivative, there is a correspondence between them.
Although the {\GL} derivative is a discrete fractional derivative and the {\RL} derivative is continuous, it was shown that both definitions are equivalent in the continuous limit \cite{kilbas2006theory, samko1993fractional}.
The {\Cp} fractional derivative can be obtained from the {\RL} fractional derivative by accounting for the initial conditions of a function at the expansion point.
The account of the initial conditions in the Caputo definition leads to a convergent form of the fractional derivative at the expansion point, in contrast to both {\GL} and {\RL} derivatives, which makes it especially suitable for physical applications \cite{herrmann2014fractional}.
In this article, we derive the exact analytical formula that casts the Gr\"{u}nwald-Letnikov fractional derivative into an infinite sum of integer order derivatives.
By representing the fractional derivative as an infinite series of integer order derivatives, we find a unified description of {\RL}, {\Cp}, and {\GL} fractional derivatives.
The only difference in our expansion for the {\RL} or {\GL} derivative and the Caputo derivative is in the lower limit of the summation index.

We examine convergence of the {\GL} fractional derivative, represented by an infinite series of integer derivatives, by truncating the infinite series and retaining only the first few terms. We find that functions normally characterized by Taylor series with a finite radius of convergence have an infinite radius of convergence in the integer derivative expansion. For physically relevant functions, such as hyperbolic tangent and secant, we show that by retaining only the first few terms in the infinite series the proposed formula efficiently approximates the fractional derivative, establishing a firm ground for its use in numerically solving fractional differential equations. Moreover, we show that for functions represented by Taylor series with an infinite radius of convergence, the truncation error is inversely proportional to the number of terms kept in the expansion. Specifically, an integer derivative expansion of $\sin(x)$ with 2 terms achieves an average $1\%$ error, and with a total of $10$ terms, the error decreases down to $0.01\%$.

Finally, we use the truncated integer derivative series to solve linear fractional differential equations with both constant and variable coefficients. We find that the fourth-order Runge-Kutta method applied to truncated fractional differential equations produces numerical solutions which rapidly converge to the exact analytical results, given by the Mittag-Leffler and generalized Fox-Wright special functions \cite{kilbas2006theory}. Approximating the fractional derivative as an integer derivative series with the first $3$ terms generates around $1\%$ error for the constant coefficient differential equation, and $10\%$ error for the 
differential equation with variable coefficients. Thus, we show that the truncated expansion provides a robust numerical scheme for solving linear fractional differential equations, such as the fractional {\Sch} and fractional diffusion equations \cite{herrmann2014fractional}.

\section{Expressing {\GL} fractional derivative as integer derivative series}
In this section we derive the infinite integer derivative expansion for the {\GL} fractional derivative. For simplicity we only consider left-sided derivatives of order $q$, with $q \in \mathds{C}$, subject to constraint of $\operatorname{Re}(q) > 0$.

We adopt the following definition of the Gr\"{u}nwald-Letnikov derivative \cite{kilbas2006theory}:
\begin{eqnarray}\label{gl_def1}
^\mathrm{GL}\mathbf{D}^q f(x) = \lim_{h \to 0 \atop N \to \infty } \frac{1}{h^q} \sum_{j=0}^{N-1} (-1)^j \binom{q}{j} f(x-jh),
\end{eqnarray}
where $N$ is the number of gridpoints and $h$ is the grid spacing defined as $h \equiv x/N$. The infinitesimal step $h$ is a constant until we perform the continuous limit. The generalized binomial coefficient $\binom{q}{j}$ valid for non-integer $q$ is defined as \cite{kilbas2006theory, samko1993fractional},
\begin{eqnarray}
\binom{q}{j} \equiv \frac{\Gamma(q+1)}{\Gamma(j+1)\hspace{1mm} \Gamma(q-j+1)} = \frac{(-1)^{j-1} \hspace{1mm} q \hspace{1mm} \Gamma(j-q)}{\Gamma(1-q) \hspace{1mm} \Gamma(j+1)},
\end{eqnarray}
where $\Gamma(z)$ is the Euler gamma function. We note that the function $f(x-mh)$ can be expressed in terms of the function $f(x)$ via the
finite displacement, or shift, operator \cite{landau1958quantum},
\begin{eqnarray}
f(x-jh) = \mathcal{D}_{jh} \left[f(x)\right] = \left(1-h\frac{d}{dx}\right)^j f(x),
\end{eqnarray}
which can be verified directly via, e.g., the finite difference method. If we make the substitution, the {\GL} derivative becomes,
\begin{IEEEeqnarray}{l}
^{\mathrm{GL}} \textbf{D}^q f(x) =
\lim_{h \to 0 \atop N \to \infty } \frac{1}{h^q} \sum_{j=0}^{N-1} (-1)^j \binom{q}{j} \Big(1-h\frac{d}{dx}\Big)^j f(x) \\ = 
\lim_{h \to 0 \atop N \to \infty }
\sum^{N-1}_{k = 0} 
\sum_{j = k}^{N-1} 
h^{k-q} \hspace{1mm}
(-1)^{j-k} \hspace{1mm}
\binom{q}{j} \hspace{1mm}
\binom{j}{k} \hspace{1mm}
\frac{d^{k}}{d x^{k}} \hspace{1mm}
f(x),
\label{sum-w-displacement-op}
\end{IEEEeqnarray}
where we applied the Newton binomial formula to the displacement operator and exchanged the order of the summation. Now we notice that we can perform the summation of the inner series,
\begin{IEEEeqnarray}{l}
 \sum_{j=k}^{N-1}  
(-1)^{j-k}  \left(
\begin{array}{c}
 j \\
 k
\end{array}
\right) \left(
\begin{array}{c}
 q \\
 j
\end{array}
\right) 
=
\frac{(-1)^{N-k + 1} (N-k) }{q-k}
 \left(
\begin{array}{c}
 N \\
 k
\end{array}
\right) \left(
\begin{array}{c}
 q \\
 N
\end{array}
\right).
\end{IEEEeqnarray}
We point out that below we treat the integer case of $q \in \mathds{N}$ separately (see Eq.(\ref{integercase1})). Thus, we cast the {\GL} derivative into an infinite sum of integer order derivatives,
\begin{eqnarray}
^{\mathrm{GL}} \textbf{D}^q f(x)  = 
\lim_{h \to 0 \atop N \to \infty }\sum_{k = 0}^{N-1} \frac{(-1)^{N-k+1} \hspace{1mm} (N-k) \hspace{1mm} h^{k-q}}{q-k} \hspace{1mm} \binom{N}{k} \hspace{1mm} \binom{q}{N} \hspace{1mm} \frac{d^k}{dx^k} \hspace{1mm} f(x).
\end{eqnarray}
To perform the limits $h \to 0$ and $N \to \infty$, we explore the weight function of the integer derivative, which we define as,
\begin{eqnarray}\label{weightfunc1}
W(q,k,N) = 
\frac{(-1)^{N-k+1} \hspace{1mm} (N-k) \hspace{1mm} h^{k-q}}{q-k} \binom{N}{k} \binom{q}{N},
\end{eqnarray}
and expand it for $N\gg{1}$ in a series in $1/N$. We obtain,
\begin{equation}
W(q,k,N) = 
(-1)^{N-k} \hspace{1mm} \left(\frac{1}{N}\right)^{q-k} \hspace{1mm} \frac{h^{k-q} \hspace{1mm} \sin\left[\pi(N-q)\right] \hspace{1mm} \Gamma(q+1)}{\pi \hspace{1mm} \Gamma(k+1)} \left(\frac{1}{k-q}-\frac{k+q+1}{2 \hspace{1mm} N}+\mathcal{O}\left(\frac{1}{N}\right)^2\right).
\end{equation}
The leading term in this expansion can be simplified to 
\begin{equation}
\lim_{h\to 0 \atop N\to \infty}
W(q,k,N) = 
\lim_{h\to 0 \atop N\to \infty}
\frac{\hspace{1mm} \sin\left[\pi(q-k)\right] \hspace{1mm} \Gamma (q+1)}{\pi \hspace{1mm} (q-k) \hspace{1mm} \Gamma(k+1)} \hspace{1mm} (hN)^{k-q}  = 
\frac{\hspace{1mm} \sin\left[\pi(q-k)\right] }{\pi \hspace{1mm} (q-k)}
\frac{\Gamma (q+1)}{\Gamma (k+1)}
 \hspace{1mm} x^{k-q}.
\end{equation}
Thus, we have performed the expansion of the {\GL} fractional derivative in terms of integer order derivatives in the limiting case of infinitesimally small grid size $h \to 0$. Lastly, we obtain our final series,
\begin{equation}\label{main-result}
^{\mathrm{GL}} \textbf{D}^q f(x) = \sum^{\infty}_{k=0} 
\frac{\hspace{1mm} \sin\left[\pi(q-k)\right] }{\pi \hspace{1mm} (q-k)}
\frac{\Gamma (q+1)}{\Gamma (k+1)}
 \hspace{1mm} x^{k-q} 
\frac{d^{k}}{d x^{k}}
f(x) = 
\sum_{k=0}^{\infty} \binom{q}{k} \frac{x^{k-q}}{\Gamma (k-q+1)} \frac{d^k}{dx^k} f(x).
\end{equation}
We note that for integer $q\in{\mathds{N}}$,
the expansion reduces to a single term due to the delta-function behavior of ${\rm sinc}[\pi(q-k)]$. Indeed for integer $q=n$ we have, 
\begin{eqnarray}\label{integercase1}
\left.\frac{\hspace{1mm} \sin\left[\pi(q-k)\right] }{\pi \hspace{1mm} (q-k)} \right|_{q=n\in{\mathds{N}}}= \delta_{n,k}.
\end{eqnarray}
Thus the infinite series of integer order derivatives reduces to a single derivative of $n^{\rm th}$ order, 
\begin{eqnarray}
\left.
^{\mathrm{GL}} \textbf{D}^q f(x) 
\right|_{q=n\in{\mathds{N}}} = 
\frac{d^{n}}{d x^{n}}
f(x)
.
\end{eqnarray}

\section{Unified description of fractional derivatives in terms of the infinite series of integer order derivatives}

In this section we establish a connection between {\RL}, {\Cp}, and {\GL} derivatives. The {\RL} fractional derivative is defined as a convolution integral, 
\begin{equation}\label{rl_def1}
^{\mathrm{RL}} \textbf{D}^q f(x) = \frac{1}{\Gamma (n-q)} \frac{d^n}{dx^n} \int_a^x 
{dt}\;
(x-t)^{n-q-1} f(t),
\end{equation}
where $n$ is the ceiling of the fractional order, $n=\ceil{q}$, given in terms of the integer part $[q]$ of $q$ as $\ceil{q} = [q]+1$. By rewriting the fractional {\RL} derivative of order $q$ as a sequential operation of an integer derivative of order $[q]+1$ and a fractional {\RL} derivative of order $q-\ceil{q}$, with a subsequent term-by-term fractional differentiation,
%, followed by fractional {\RL} derivative of an order $q- \ceil{q}$, with a subsequent term-by-term fractional differentiation,
one obtains an integer derivative expansion for the {\RL} fractional derivative  \cite{samko1993fractional},
\begin{eqnarray}\label{rl_result1}
 ^{\mathrm{RL}} \textbf{D}^q f(x) =
\sum^{\infty}_{k=0}
\binom{q}{k}
\frac{(x-a)^{k-q}}{\Gamma(k-q+1)} f^{(k)}(x),
\end{eqnarray}
where $a$ is the base point of the {\RL}
derivative. By choosing a zero base point $a=0$ and comparing our formula Eq.(\ref{main-result}) for the {\GL} derivative with the expansion derived for the {\RL} derivative Eq.(\ref{rl_result1}), we conclude that the {\GL} fractional derivative given by Eq.(\ref{gl_def1}) 
and the {\RL} fractional derivative given by Eq.(\ref{rl_def1})  are not only equivalent in the continuous limit but also lead to the very same infinite expansion of integer order derivatives.

In order to obtain a unified description for all three fractional derivatives, we introduce the
{\Cp} fractional derivative, defined according to \cite{kilbas2006theory},
\begin{equation}
^{\mathrm{C}} \textbf{D}^q f(x) = \frac{1}{\Gamma (n-q)} \int_a^x {dt}\; (x-t)^{n-q-1} \frac{d^n f(t)}{dt^n}.
\end{equation}
First, we refer to the connection between {\Cp} and {\RL} fractional derivatives in \cite{kilbas2006theory},
\begin{IEEEeqnarray}{l}
^{\mathrm{C}} \textbf{D}^{q} f(x) = {}^{\mathrm{RL}} \textbf{D}^q
\left[
f(x)- 
\sum^{n-1}_{k=0}
\frac{f^{(k)}(a)}{k!}(x-a)^{k}
\right]
,
\end{IEEEeqnarray}
where integer order derivatives are evaluated at the base point $a$, i.e.
\begin{eqnarray}
f^{(k)}(a)\equiv{}\left.\frac{d^{k}f(x)}{dx^{k}}\right|_{x=a}.
\end{eqnarray}
By applying the infinite expansion in Eq.(\ref{main-result}), we obtain,
\begin{equation}\label{cap_result1}
{}^{\mathrm{C}}{\mathbf{D}}^{q}f(x)\equiv{}\lim_{N\to\infty}{{}^{\mathrm{C}}\mathds{D}}_{N}^{q}f(x) =  \lim_{N\to\infty}\sum^{N}_{j=0} \binom{q}{j}
\frac{(x-a)^{j-q}}{\Gamma(j-q+1)}
\left[
f^{(j)}(x) 
-
\sum^{n-1}_{k=j}
\frac{f^{(k)}(a)}{(k-j)!}
(x-a)^{k-j}
\right].
\end{equation}
We see that Eq.(\ref{main-result}) is versatile because it bundles all three fractional derivatives into a single expansion, with a simple adjustment on the lower bound for the Caputo derivative series and a zero base point on the {\RL} and Caputo fractional derivatives. Thus, the integer derivative expansion in Eqs.(\ref{main-result}), (\ref{rl_result1}), and Eq.(\ref{cap_result1})
gives a universal formulation for all three fractional derivatives. This universality is an important consistency test for fractional calculus. Moreover, the infinite expansion in Eq.(\ref{main-result}) is particularly convenient for numerical implementation in linear FDEs, as we present in the following sections.

\section{Truncation, Error, and Radius of Convergence}
%The formula Eq. (\ref{main-result}) gives an infinite series of integer derivative. 

In the previous section we obtained the unified description of {\GL}, {\RL}, and {\Cp} fractional derivatives in terms of an infinite series of integer order derivatives. Even though the infinite expansion of the {\RL} fractional derivative was derived previously \cite{samko1993fractional}, the numerical applications of the result in Eq.(\ref{main-result}) and Eq.(\ref{rl_def1}), which necessarily rely on the truncation of the infinite series, were missing. The goal of this section is to 
truncate the infinite series given by Eq.(\ref{main-result}) and calculate the residual truncation error for several physically relevant functions. To determine the error introduced by truncating the series, we perform multiple case studies in which we consider functions with both an infinite radius of convergence  of the Taylor series, such as plane and standing waves, Gaussian function, as well as functions with a finite radius of convergence, e.g., hyperbolic secant (hyperbolic tangent) which describe bright (dark) soliton propagation. Moreover, we evaluate the minimal number of terms kept in the infinite series which correspond to a given level of accuracy.
%the number of terms that should be retained to guarantee a desired error.
In particular, we choose the Caputo fractional derivative of the order  $q = {1}/{2}$. 
We calculate relative error by
\begin{equation}\label{deferror1}
\epsilon(x) = \frac{a(x)-b(x)}{\frac{1}{2}(|a(x)|+|b(x)|)},
\end{equation}
where $a(x) ={}^{\mathrm{C}} \textbf{D}^{q} f(x)$ is the infinite series given by 
Eq.(\ref{cap_result1}) and $b(x) = {}^{\mathrm{C}} {\mathds{D}}_{N}^{q} f(x)$ is the truncated series, where  $q$ is the order of the fractional differential operator, and $N$ is the number of terms in the truncated series Eq.(\ref{cap_result1}). We observe spikes in the log-error, $\log(|\epsilon(x)|)$, either in case of real-valued roots of the fractional derivative $a(x)$ or its approximation $b(x)$, or in the case of a match between the fractional derivative and its approximation. Yet another discontinuity in the error arises if the fractional derivative and its approximation are of equal magnitude but opposite in sign.

\begin{figure}
\centering
\begin{subfigure}[b]{0.45\textwidth}
  \centering
  \captionsetup{justification=centering}
  \vspace{2mm}
\includegraphics[width=\textwidth]{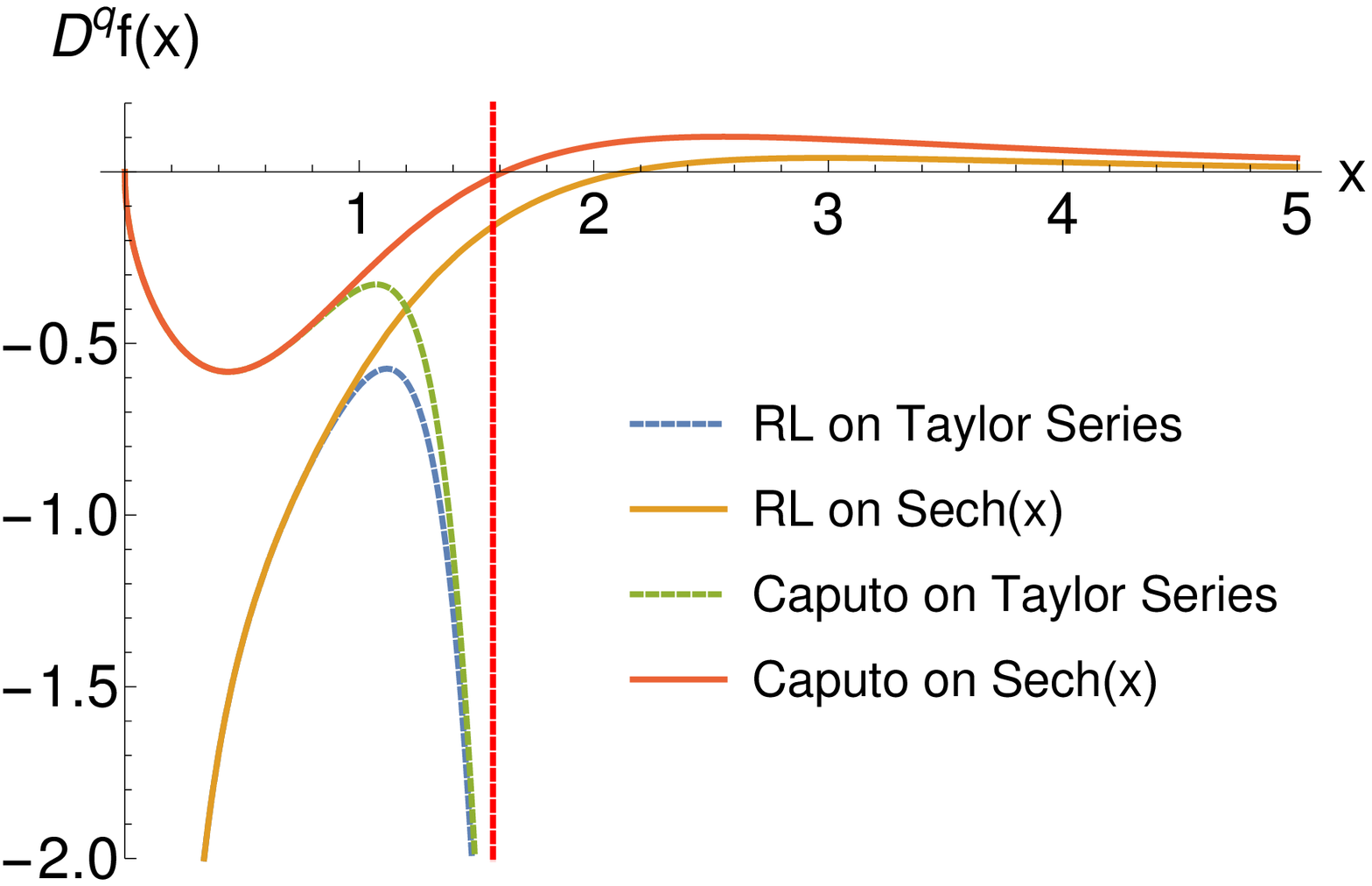}
 \caption{•}
 \label{sech-and-its-taylor-series}
\end{subfigure}%
\begin{subfigure}[b]{0.55\textwidth}
  \centering
  \captionsetup{justification=centering}
  \vspace{4mm}
\includegraphics[width=\textwidth]{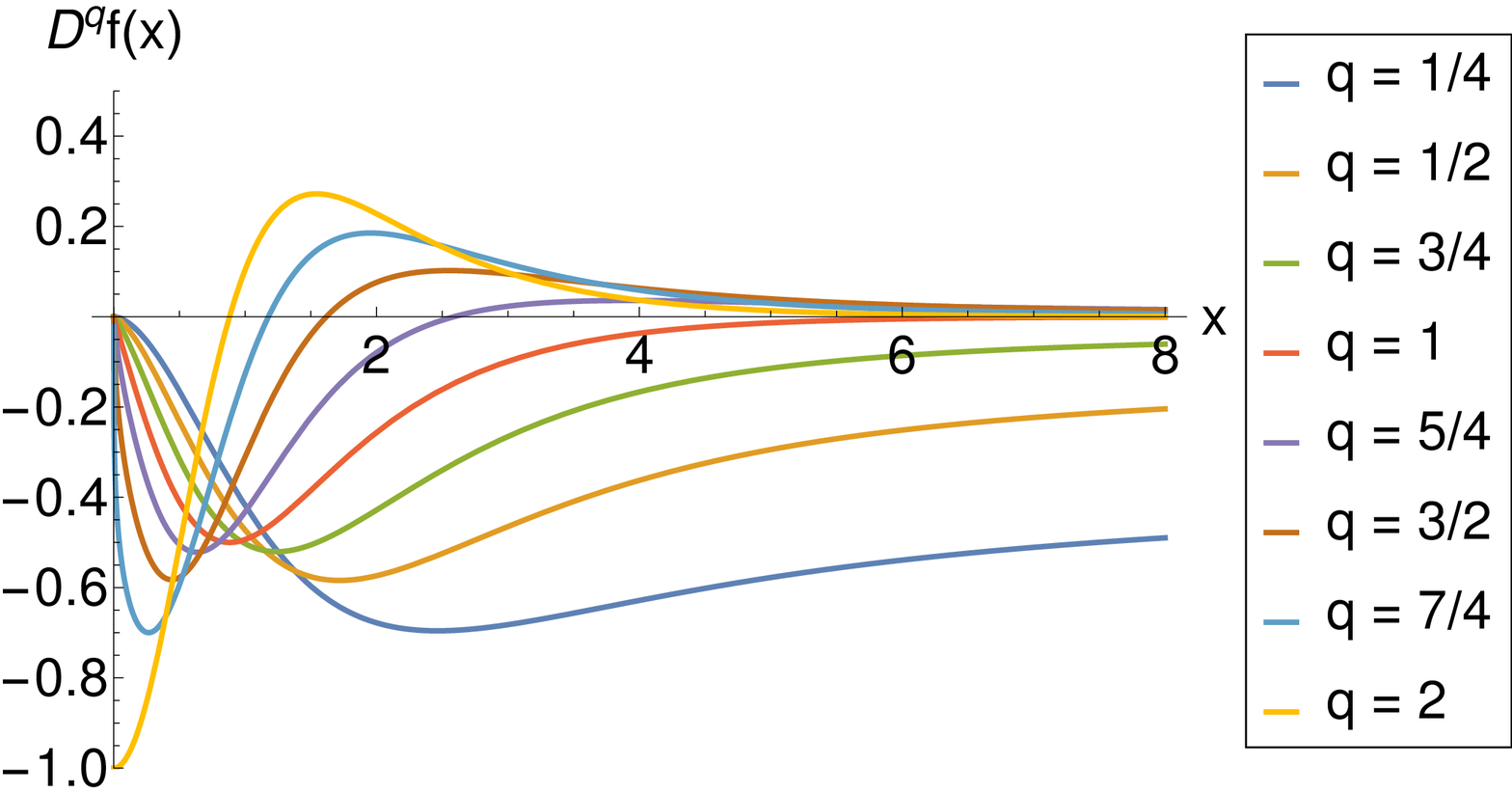}
 \caption{•}
 \label{sec-many-q}
\end{subfigure}%
\\
\begin{subfigure}[b]{0.45\textwidth}
  \centering
  \captionsetup{justification=centering}
\includegraphics[width=\textwidth]{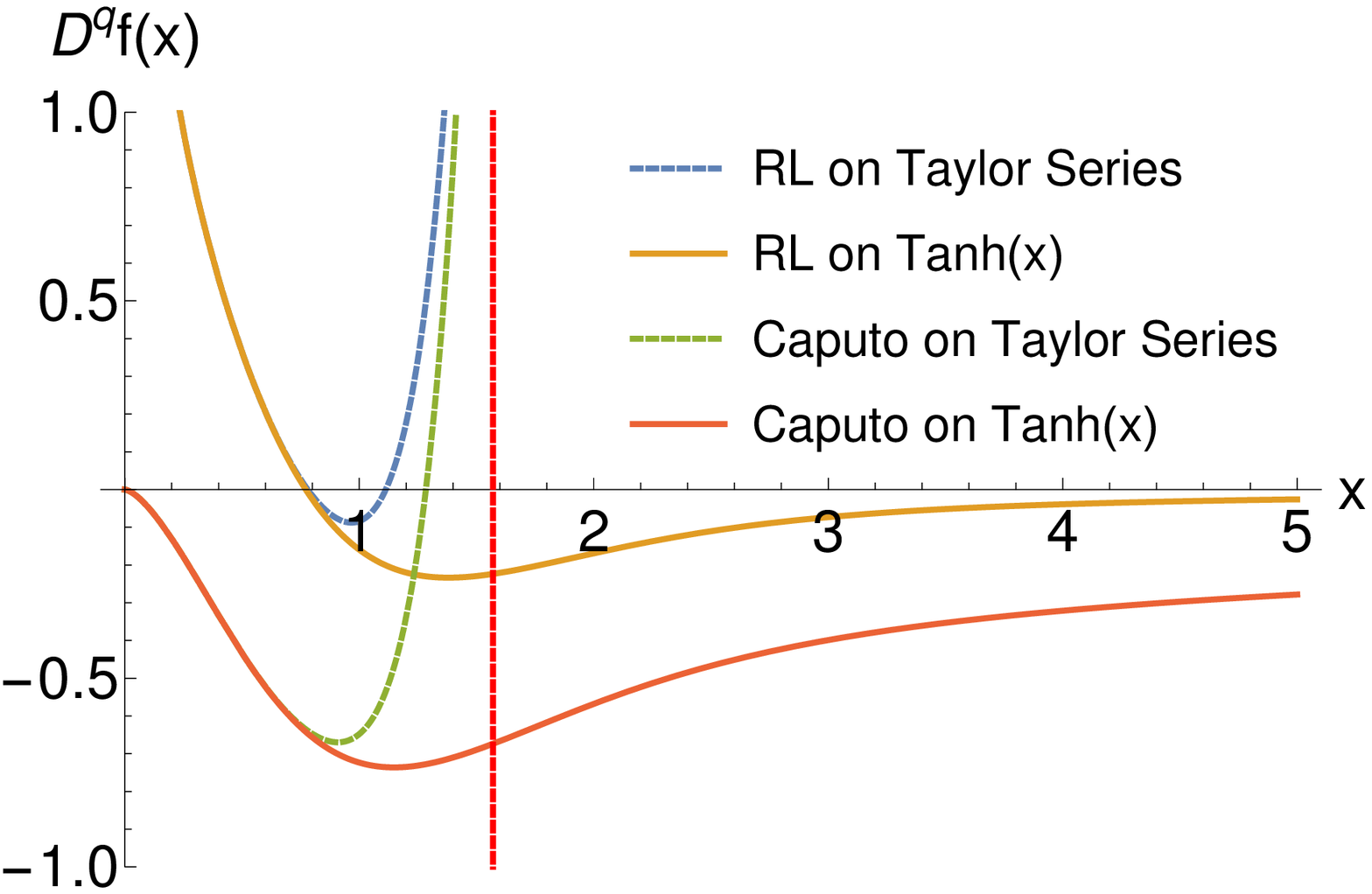}
 \caption{•}
 \label{tanh-and-its-taylor-series}
\end{subfigure}%
\begin{subfigure}[b]{0.55\textwidth}
  \centering
  \captionsetup{justification=centering}
\includegraphics[width=\textwidth]{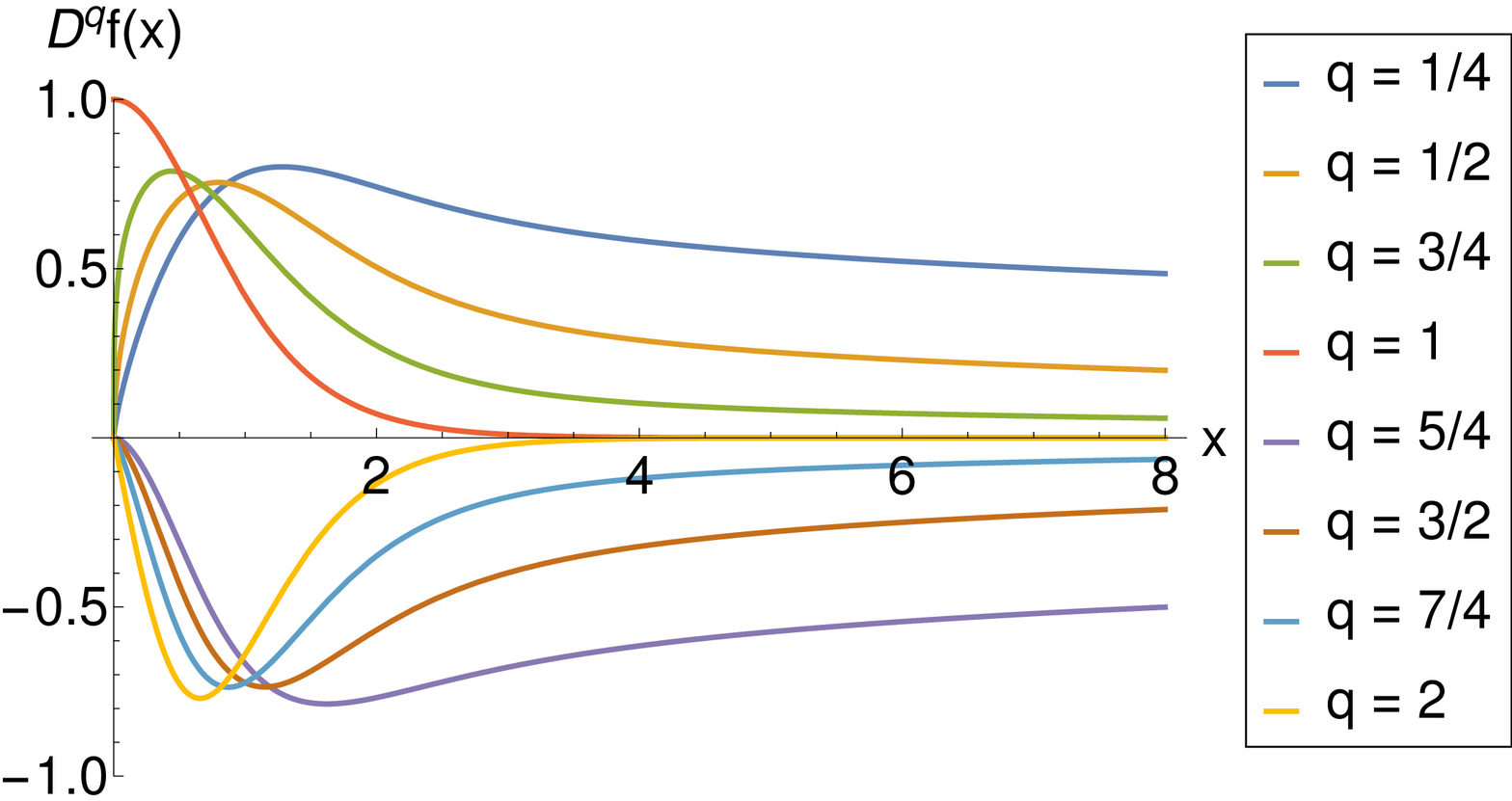}
 \caption{•}
 \label{tanh-many-q}
\end{subfigure}%
\caption{
Fractional derivative of hyperbolic secant and tangent functions.
(a) {\RL} fractional derivative (blue curve) and {\Cp} fractional derivative (green curve) of order $q=3/2$ evaluated via Taylor expansion of hyperbolic secant function are divergent at $R=\pi/2$ due to finite radius of convergence of the Taylor series. However, our representation of {\RL} 
fractional derivative (orange curve) and  {\Cp} fractional derivative (red curve) in terms of an infinite series of integer derivatives of the original function given by Eq.(\ref{main-result}) does not rely on properties of its Taylor series and, thus, leads to an infinite radius of convergence.
(b) Infinite series representation of the {\Cp} derivative of hyperbolic secant
is convergent for a whole range of fractional orders $1/4\leq{q}\leq{2}$ (shown in the legend) beyond the radius of convergence of its Taylor series $R=\pi/2$.
(c) Same as (a) but for hyperbolic tangent.
(d) Same as (b) but for hyperbolic tangent.
}
%Integer derivative series applied to $\sech(x)$ and Taylor series of $\sech(x)$ for $q = \frac{3}{2}$. The integer derivative expansion contains 4 terms. We see that when we expand $\sech(x)$ in a Taylor series, its fractional derivative through the integer derivative expansion inherits the finite radius of convergence at $\pi/2$ (marked in red dash).}
\label{finite-radius-of-converge}
\end{figure}

\begin{figure}[h!]
\centering
\hspace{-10mm}
\begin{subfigure}[b]{0.5\textwidth}
  \centering
  \captionsetup{justification=centering}
\includegraphics[width=\textwidth]{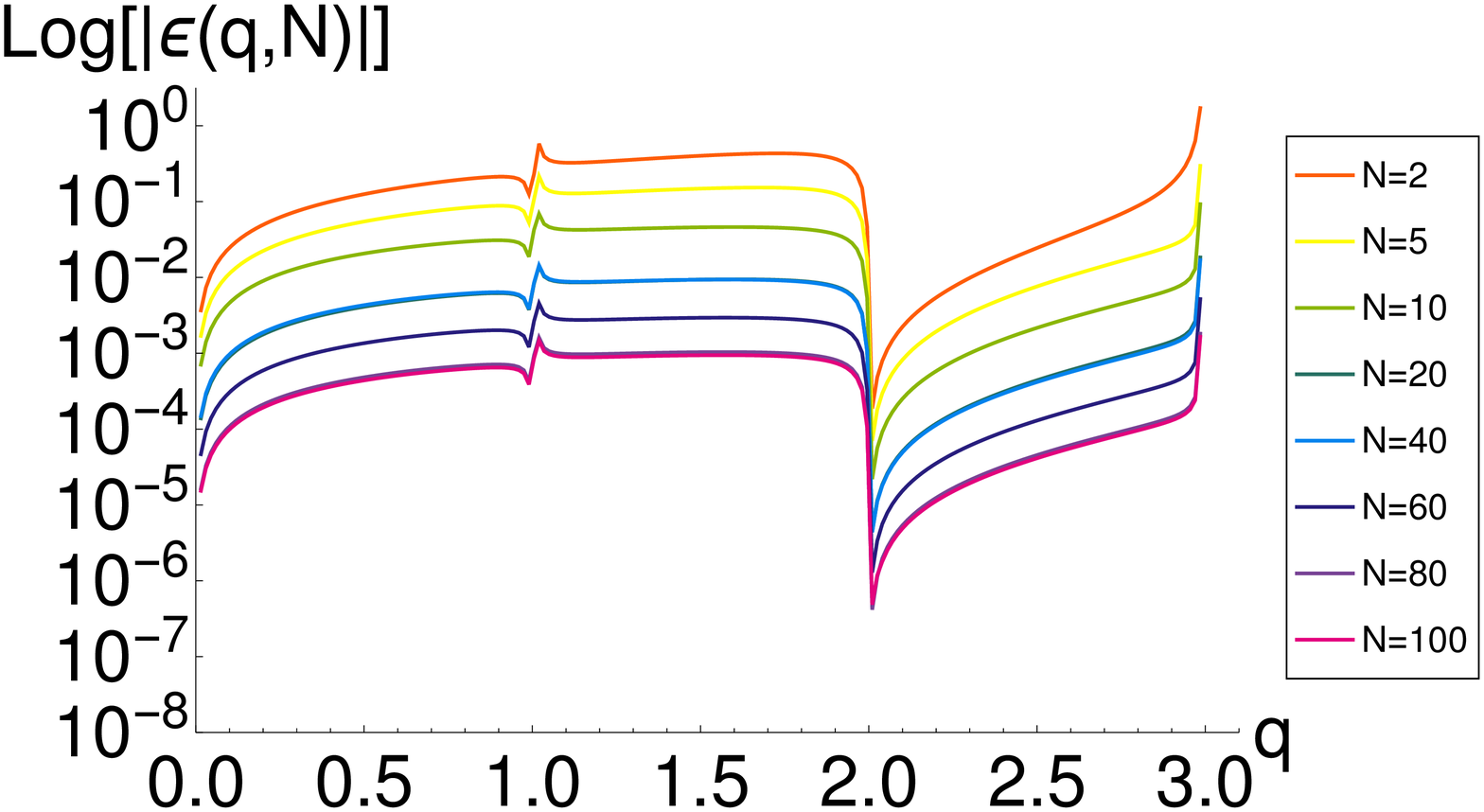}
 \caption{•}
\end{subfigure}%
\begin{subfigure}[b]{0.5\textwidth}
  \centering
  \captionsetup{justification=centering}
\includegraphics[width=\textwidth]{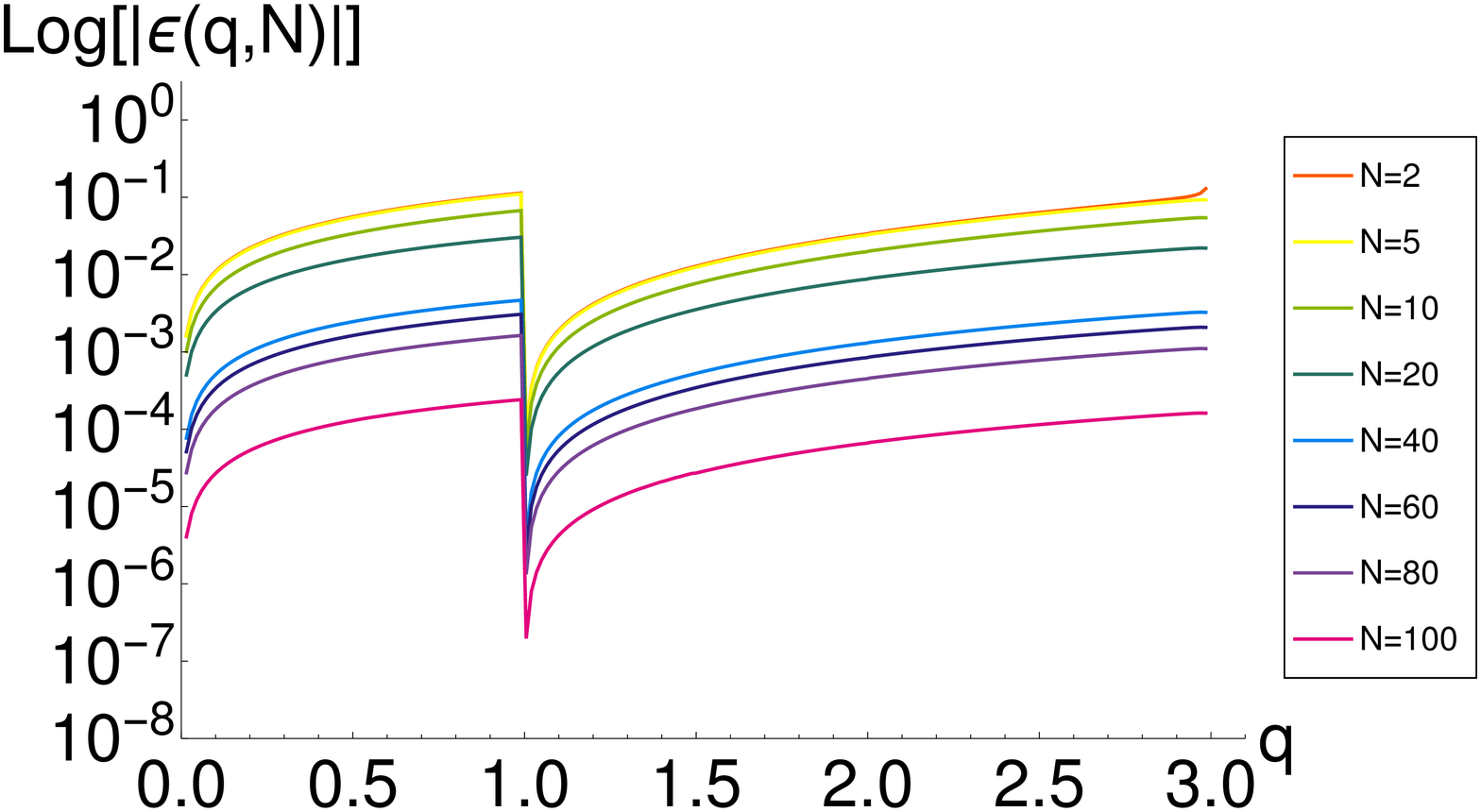}
 \caption{•}
\end{subfigure}%
%\vspace{-40mm}
\caption{
Log-linear plot of the truncation error in the fractional derivative of (a) $f(x)=\sech(x)$ 
and (b) $f(x)=\tanh(x)$ as a function of fractional order $q$, and number of terms $N$ kept in the infinite expansion Eq.(\ref{main-result}), shown in the legend.}
\label{log-linear-error-sech-tanh}
\end{figure}

To approximate the fractional derivative of hyperbolic secant to within $10\%$, we need to keep only the first three terms, as can be seen in Fig.(\ref{log-linear-error-sech-tanh}). We note that
a traditional approach in the evaluation of a fractional derivative of hyperbolic secant
(tangent) relies on the Taylor series expansion, which  diverges at $R=\pi/2$ due to a pole in the complex plane \cite{brown2009complex}. The divergence of the Taylor series results in the divergence of {\RL} and {\Cp} fractional derivatives if it is directly used in the integration process. However, the infinite series representation of the fractional derivative of $\sech(x)$ and $\tanh(x)$ given by Eq.(\ref{main-result}) is formulated in terms of integer derivatives of the original function, and does not depend on the properties of the Taylor series. Thus, the integer derivative series for the {\GL} fractional derivative of hyperbolic secant and hyperbolic tangent functions has an infinite radius of convergence, as can be seen in Fig.(\ref{finite-radius-of-converge}). The log-linear plot of truncation error in the fractional derivative of $\sech(x)$ and $\tanh(x)$ is shown in Fig.(\ref{log-linear-error-sech-tanh}).
 
\begin{figure}[h!]\label{comberror}
\centering
\begin{subfigure}[b]{0.5\textwidth}
  \centering
  \captionsetup{justification=centering}
\includegraphics[width=\textwidth]{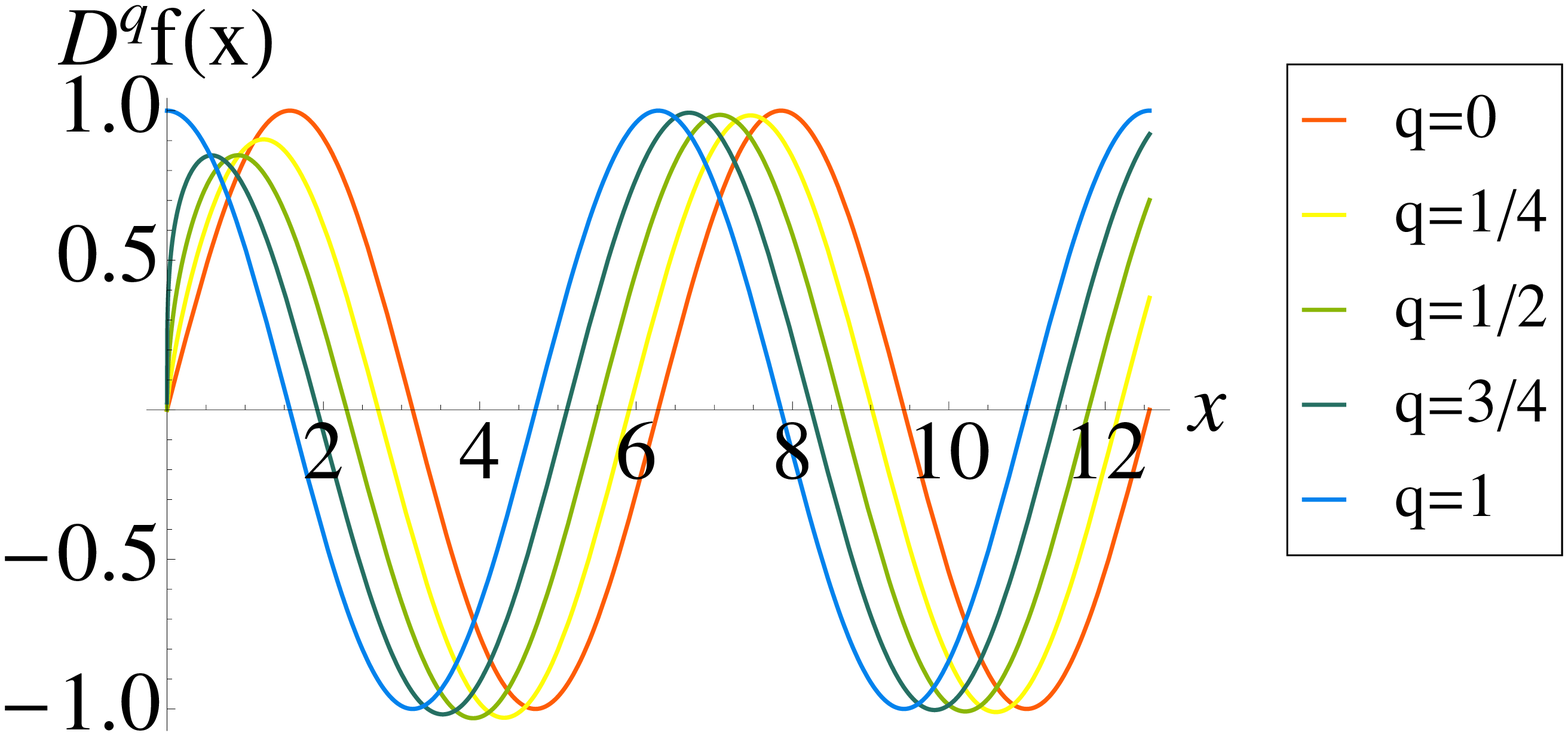}
 \caption{•}
\end{subfigure}%
\begin{subfigure}[b]{0.5\textwidth}
  \centering
  \captionsetup{justification=centering}
\includegraphics[width=\textwidth]{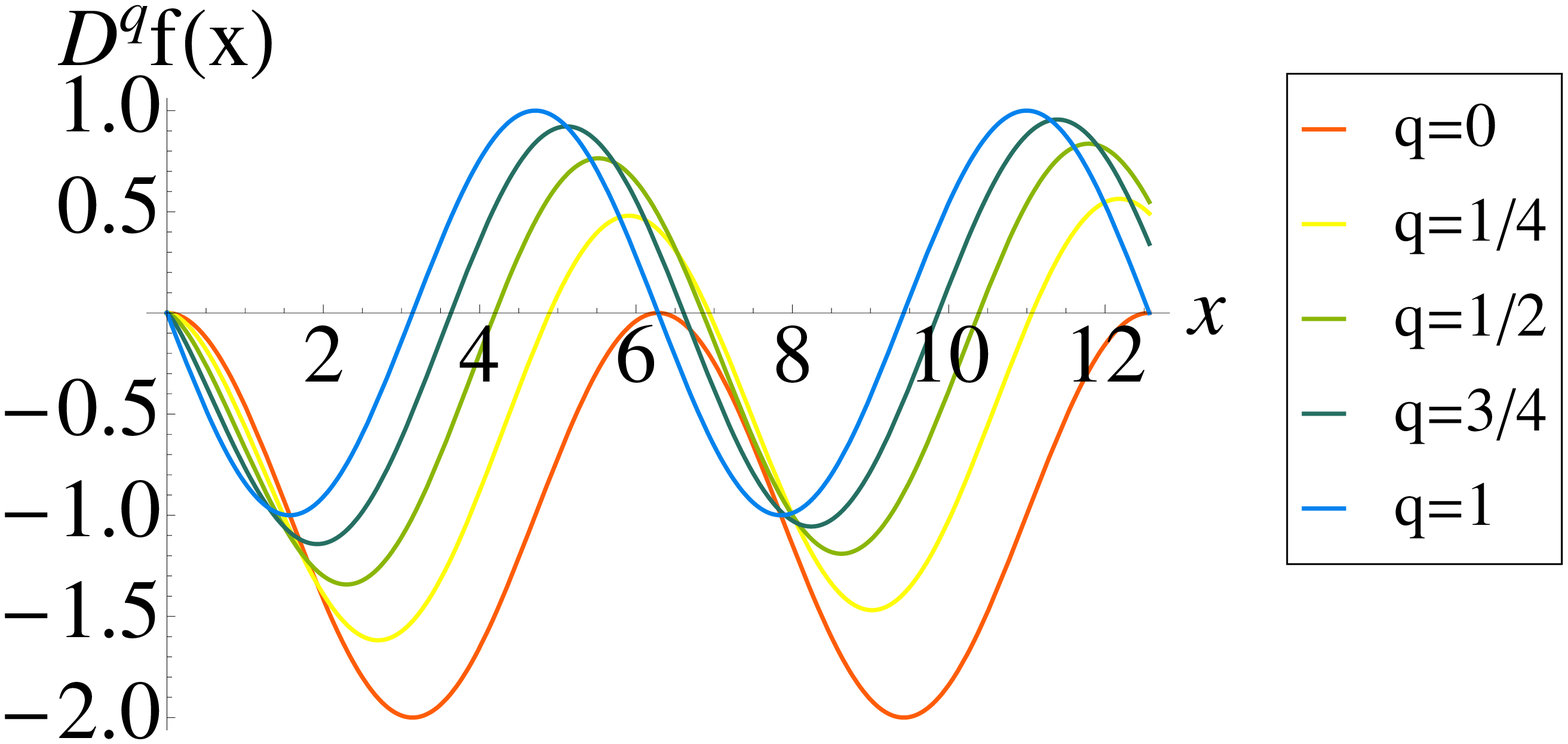}
 \caption{•}
\end{subfigure}%
\\
\begin{subfigure}[c]{0.5\textwidth}
  \centering
  \captionsetup{justification=centering}
\includegraphics[width=\textwidth]{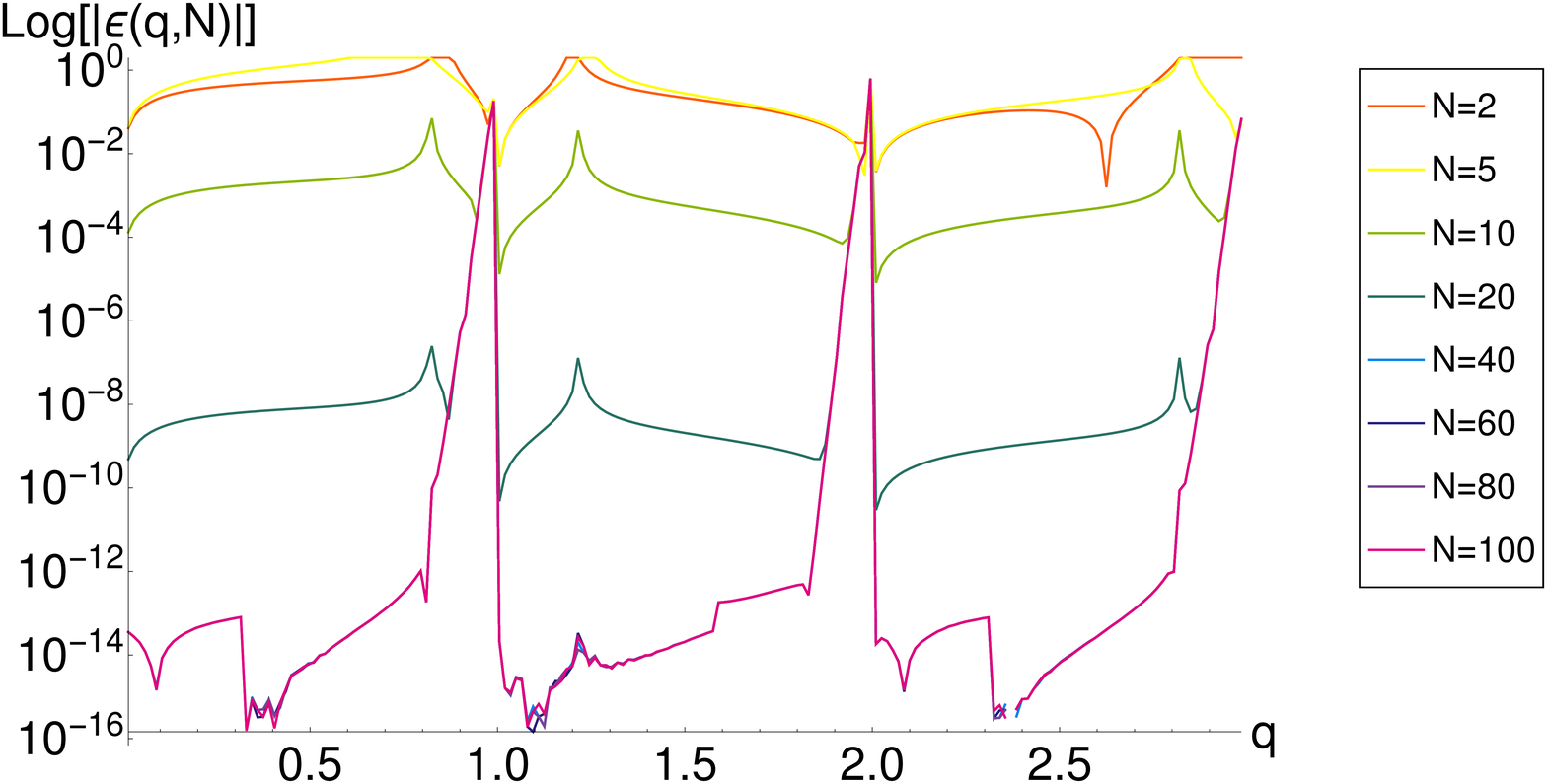}
 \caption{•}
\end{subfigure}%
\begin{subfigure}[d]{0.5\textwidth}
  \centering
  \captionsetup{justification=centering}
\includegraphics[width=\textwidth]{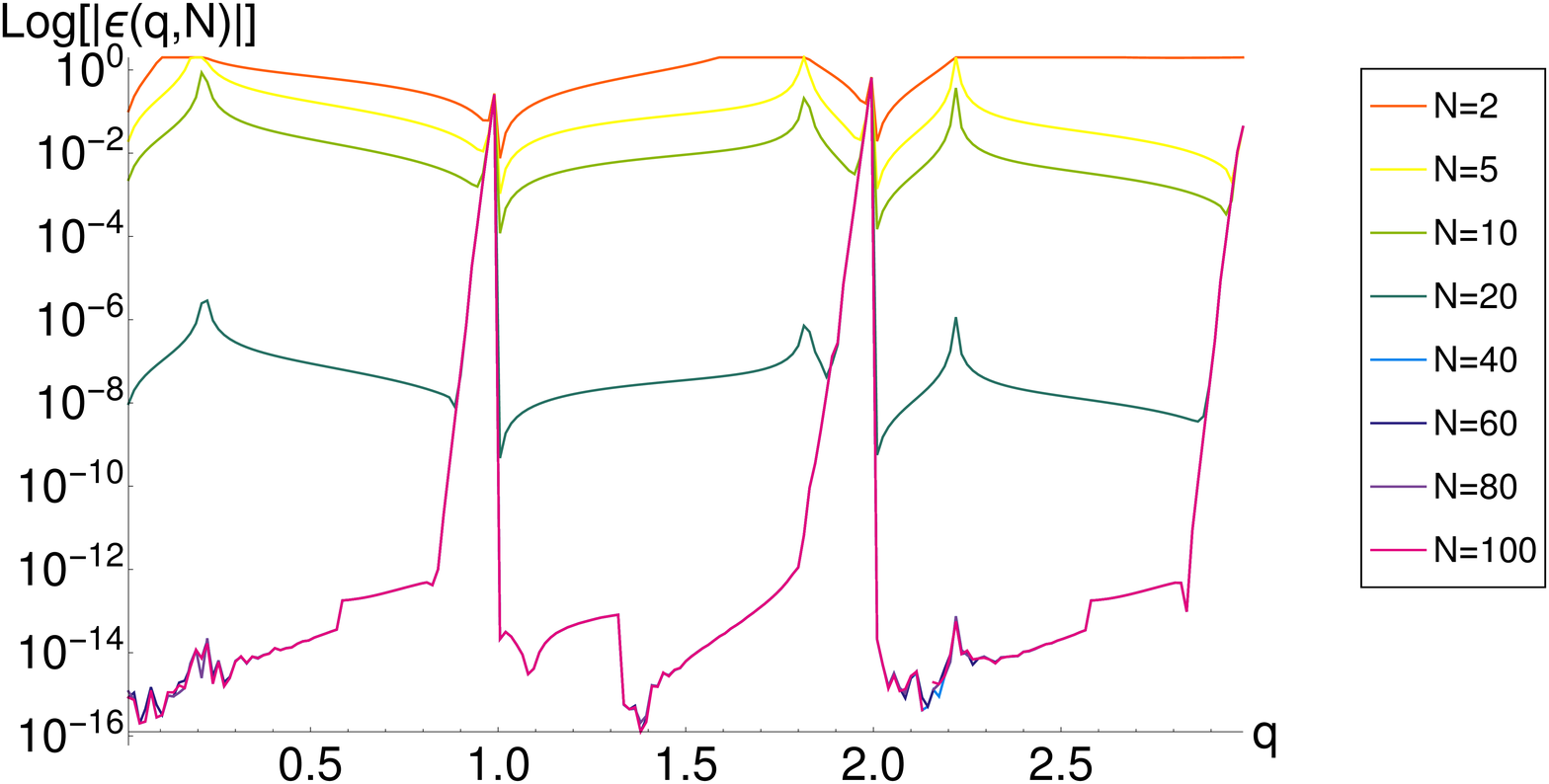}
 \caption{•}
\end{subfigure}%
\caption{
Caputo fractional derivative of (a) $\sin(x)$ and (b) $\cos(x)$ as a function of position $x$ and fractional order $q$ in the range $0\leq{q}\leq{1}$ (shown in legend).
Log-linear plot of the truncation error in the fractional derivative of (c) $\sin(x)$
and (d) $\cos(x)$ as a function of fractional order $q$, and number of terms $N$ kept in the infinite expansion Eq.(\ref{main-result}), shown in the legend.
}
\label{sin-cos-plots}
\end{figure}

For functions described by Taylor series with an infinite radius of convergence, e.g., $\sin(x)$ and $\cos(x)$, the number of terms needed to reach a given level of accuracy depends on the distance away from the base point used in the integer derivative expansion. For example, to approximate the Caputo fractional derivative on $\cos(x)$, we need to retain the first $15$ terms to reach $10\%$ accuracy in the same domain as for the fractional derivative on the hyperbolic secant, as can be seen in Fig.(\ref{sin-cos-plots}).

While for a certain class of functions the integer derivative series given by Eq.(\ref{main-result}) improves the fractional derivative approximation with every additional term, the integer derivative expansion of a fractional derivative of a  Gaussian function diverges for finite orders of $N$, as we show below in Fig.(\ref{comberror}). We note that the finite sum is convergent only in a vicinity around the origin and at infinity due to the Gaussian envelope. Indeed, our integer derivative expansion given by Eq.(\ref{main-result}) for the {\GL} fractional derivative of $\exp(-x^2)$ can be expressed in terms of Hermite polynomials $H_k(x)$, i.e.
\begin{IEEEeqnarray}{l}\label{hermite1}
^{\mathrm{GL}} {\mathds{D}}_{N}^{q} [e^{-x^2} ]  = 
\sum^{N-1}_{{k}=0} 
\frac{\sin\left[\pi(q-{k})\right] }{\pi \hspace{1mm} (q-{k})}
\frac{\Gamma (q+1)}{\Gamma ({k}+1)}
 \hspace{1mm} x^{{k}-q} 
\frac{d^{k}}{dx^{k}}e^{-x^2} = \\\nn = 
e^{-x^2}
\sum^{N-1}_{{k}=0}
\frac{\sin\left[\pi(q-{k})\right] }{\pi \hspace{1mm} (q-{k})}
\frac{\Gamma (q+1)}{\Gamma ({k}+1)}
 \hspace{1mm} x^{{k}-q} 
H_{k}(-x).
\end{IEEEeqnarray}
This sum inherits large oscillations from the Hermite polynomials for both large values of its argument $x$ and its index $k$. These oscillations result in a divergence of the integer derivative expansion, and thus, establish limits of the universality of the main result Eq.(\ref{main-result}).

%\begin{figure}[h!]\label{comberror}
%\centering
%\includegraphics[width=0.8\textwidth]%{tanh_many_q.eps}
%\caption{
%Fractional derivative of $\tanh(x)$ as a %function of position $x$ and fractional order $q$ in the range $1/4\leq{q}\leq{2}$ (shown in legend).}
%\end{figure}

%\begin{figure}[h!]\label{comberror}
%\centering
%\includegraphics[width=0.8\textwidth]{tanh_log_error_sheets.eps}
%\caption{
%Log of truncation error in the fractional derivative of hyperbolic tangent 
%as a function of fractional order $q$, position $x$, and number of terms kept in the 
%infinite expansion given by Eq.(\ref{main-result}) with $N=2$ (red sheet), $N=5$ (green sheet), $N=10$ (blue sheet), $N=20$ (orange sheet).}
%\end{figure}

%\clearpage

\begin{figure}[h!]
\centering
\begin{subfigure}[b]{0.4\textwidth}
  \centering
  \captionsetup{justification=centering}
\includegraphics[width=\textwidth]{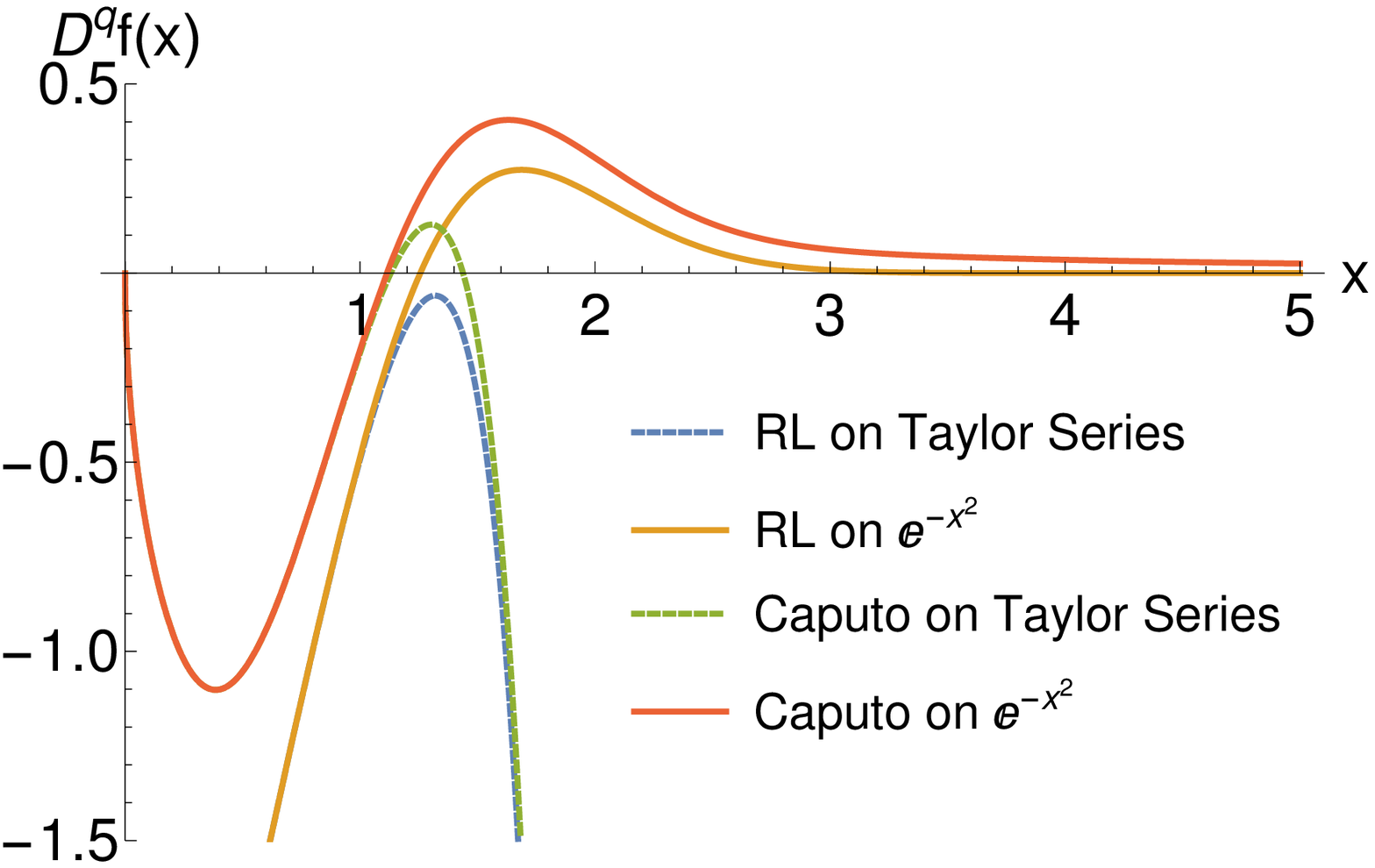}
 \caption{•}
\end{subfigure}%
\hspace{5mm}
\begin{subfigure}[b]{0.5\textwidth}
  \centering
  \captionsetup{justification=centering}
\includegraphics[width=\textwidth]{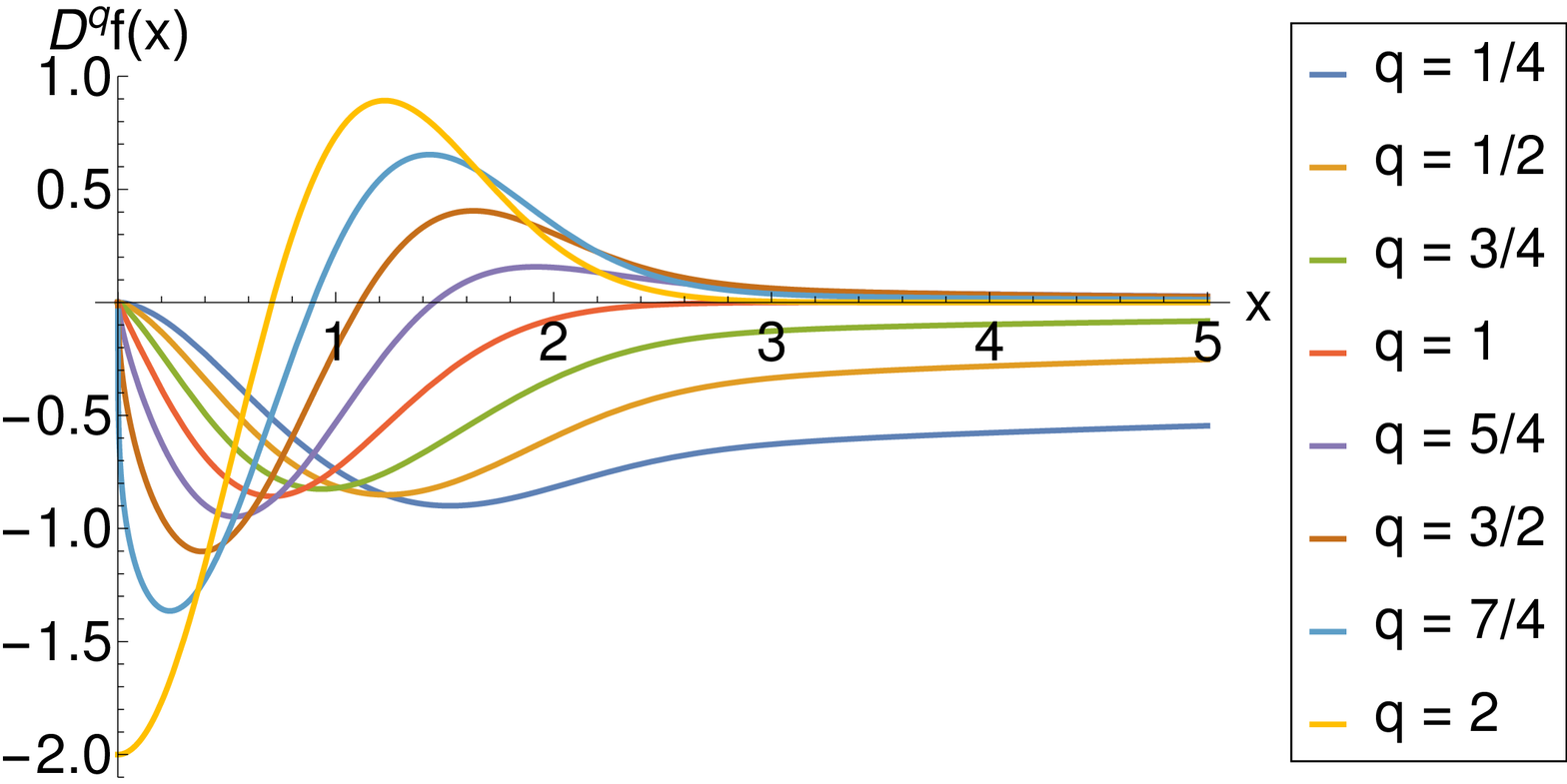}
 \caption{•}
\end{subfigure}%
\\
\centering
\begin{subfigure}[b]{0.45\textwidth}
  \centering
  \captionsetup{justification=centering}
\includegraphics[width=\textwidth]{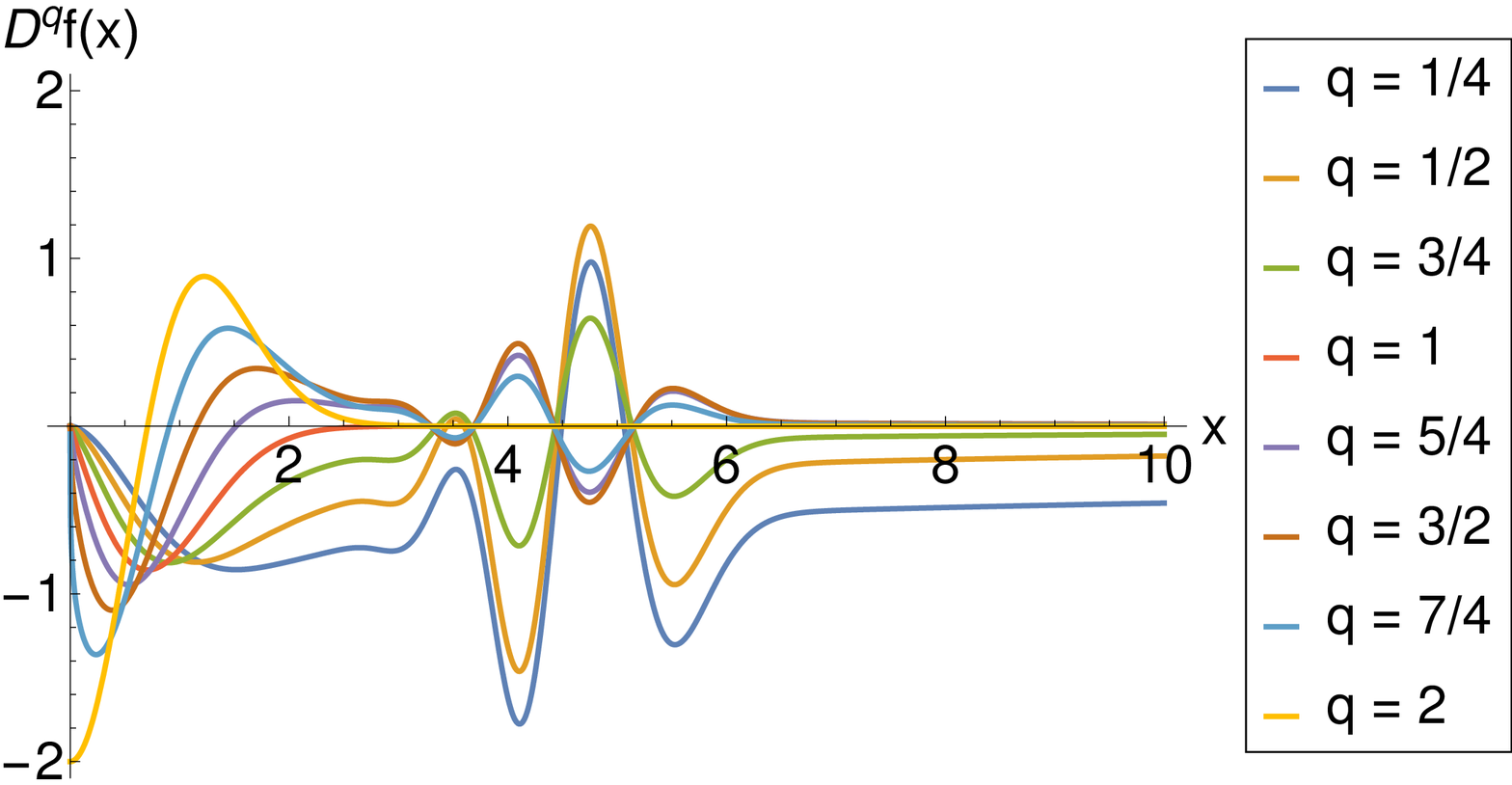}
 \caption{•}
\end{subfigure}%
\hspace{5mm}
\begin{subfigure}[b]{0.45\textwidth}
  \centering
  \captionsetup{justification=centering}
\includegraphics[width=\textwidth]{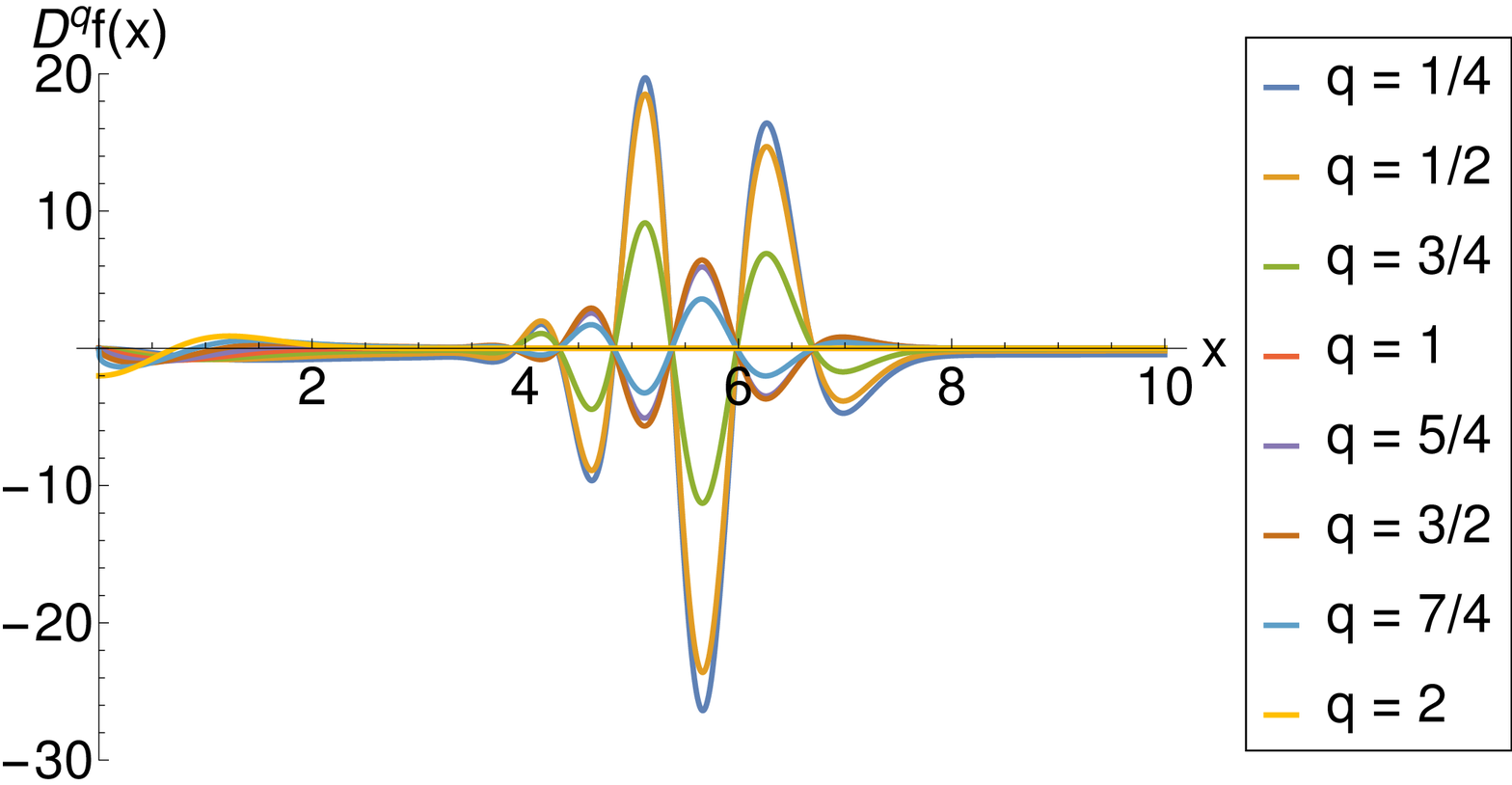}
 \caption{•}
\end{subfigure}%
\caption{
Fractional derivative of a Gaussian function.
(a) {\RL} fractional derivative (blue curve) and {\Cp} fractional derivative (green curve) of order $q=3/2$, evaluated via a Taylor expansion of a Gaussian function. As we take more terms in the Taylor expansion for $\exp(-x^2)$, the {\RL} and {\Cp} fractional derivatives converge to the orange and red curves, respectively, calculated by the integer derivative series in Eq.(\ref{main-result}).
(b)
The truncated expansion Eq.(\ref{main-result}) of the Caputo fractional derivative of a Gaussian function with only $N=3$ terms. The integer derivative series for a Gaussian function (see Eq.(\ref{hermite1})) can be written in terms of Hermite polynomials $H_{n}(x)$ which oscillate and grow factorially with $n\to\infty$.
As a consequence, the integer derivative expansion Eq.(\ref{main-result}) for a Gaussian with $N\gg{q}$ is divergent as can be directly seen in (c) which shows the truncated expansion with $N=20$ terms and (d) $N=40$ terms.}
\label{comberror}
\end{figure}

%\begin{figure}[h!]\label{comberror}
%\centering
%\includegraphics[width=0.8\textwidth]{EXACT_CAPUTO_gaussian_log_error_sheets.eps}
%\caption{
%Log of truncation error in the fractional derivative of a Gaussian 
%as a function of fractional order $q$, position $x$, and number of terms kept in the 
%infinite expansion given by Eq.(\ref{main-result}) with $N=2$ (red sheet), $N=5$ (green sheet), $N=10$ (blue sheet), $N=20$ (orange sheet).}
%\end{figure}

\clearpage

%\begin{figure}
%\centering
%\begin{subfigure}[b]{.5\textwidth}
%  \centering
%  \captionsetup{justification=centering}
% \includegraphics[width=\textwidth]{sine_approx.eps}
% \caption{•}
%  %\caption{Purple dash is approximation, solid cyan is exact.}
%  %\label{sine-approx}
%\end{subfigure}%
%\begin{subfigure}[b]{.5\textwidth}
%  \centering
%  \captionsetup{justification=centering}
%  \includegraphics[width=\textwidth]{sine_error.eps}
%  \caption{•}
%  %\caption{Solid purple is error, pink dash is 10 percent mark.}
%  %\label{sine-err}
%\end{subfigure}
%\caption{Approximation of $q=\frac{1}{2}$ Caputo derivative on $f(x) = \sin(x)$ with 15 terms. Gray dash indicates the two points where the fractional derivative and its approximation cross the $x$-axis. (a) Purple dash is approximation, solid cyan is exact. (b) Solid purple is error, pink dash is 10 percent mark.}
%\end{figure}
%
%\begin{figure}
%\centering
%\begin{subfigure}{.5\textwidth}
%  \centering
%  \includegraphics[width=0.95\linewidth]{exp_approx.eps}
%  \caption{Purple dash is approximation, solid cyan is exact.}
%  \label{exp-approx}
%\end{subfigure}%
%\begin{subfigure}{.5\textwidth}
%  \centering
%  \includegraphics[width=0.95\linewidth]{exp_error.eps}
%  \caption{Solid purple is error, pink dash is 10 percent mark.}
%  \label{exp-err}
%\end{subfigure}
%\caption{Approximation of $q=\frac{1}{2}$ Caputo derivative on $f(x) = \exp(x)$ with 15 terms.}
%\end{figure}

\section{Solving linear fractional differential equations with constant and variable coefficients using truncated series}
\begin{figure}[h!]
\centering
\begin{subfigure}[b]{0.45\textwidth}
  \centering 
  \captionsetup{justification=centering}
\includegraphics[width=\textwidth]{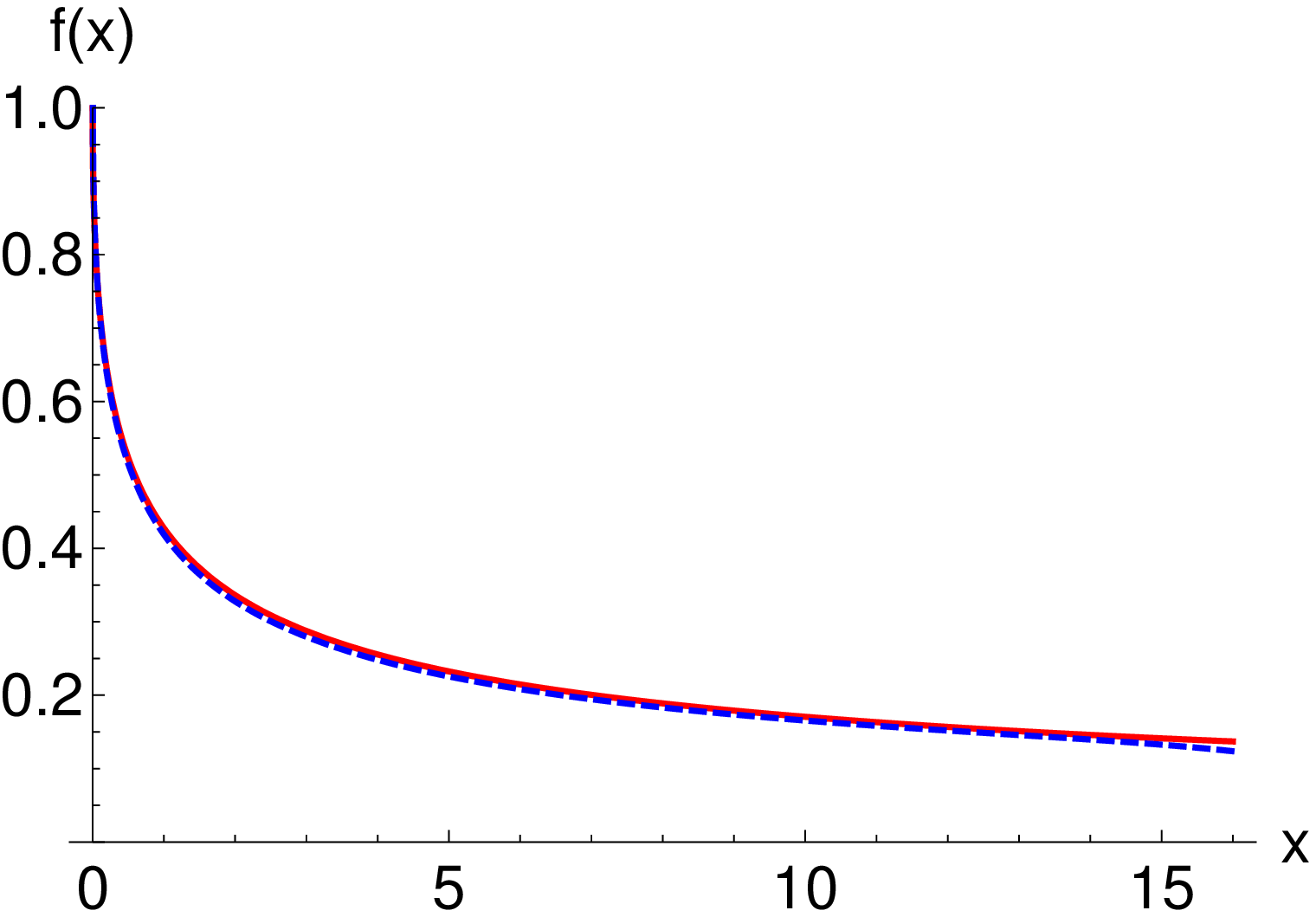}
 \caption{•}
\end{subfigure}%
\begin{subfigure}[b]{0.45\textwidth}
  \centering
  \captionsetup{justification=centering}
\includegraphics[width=\textwidth]{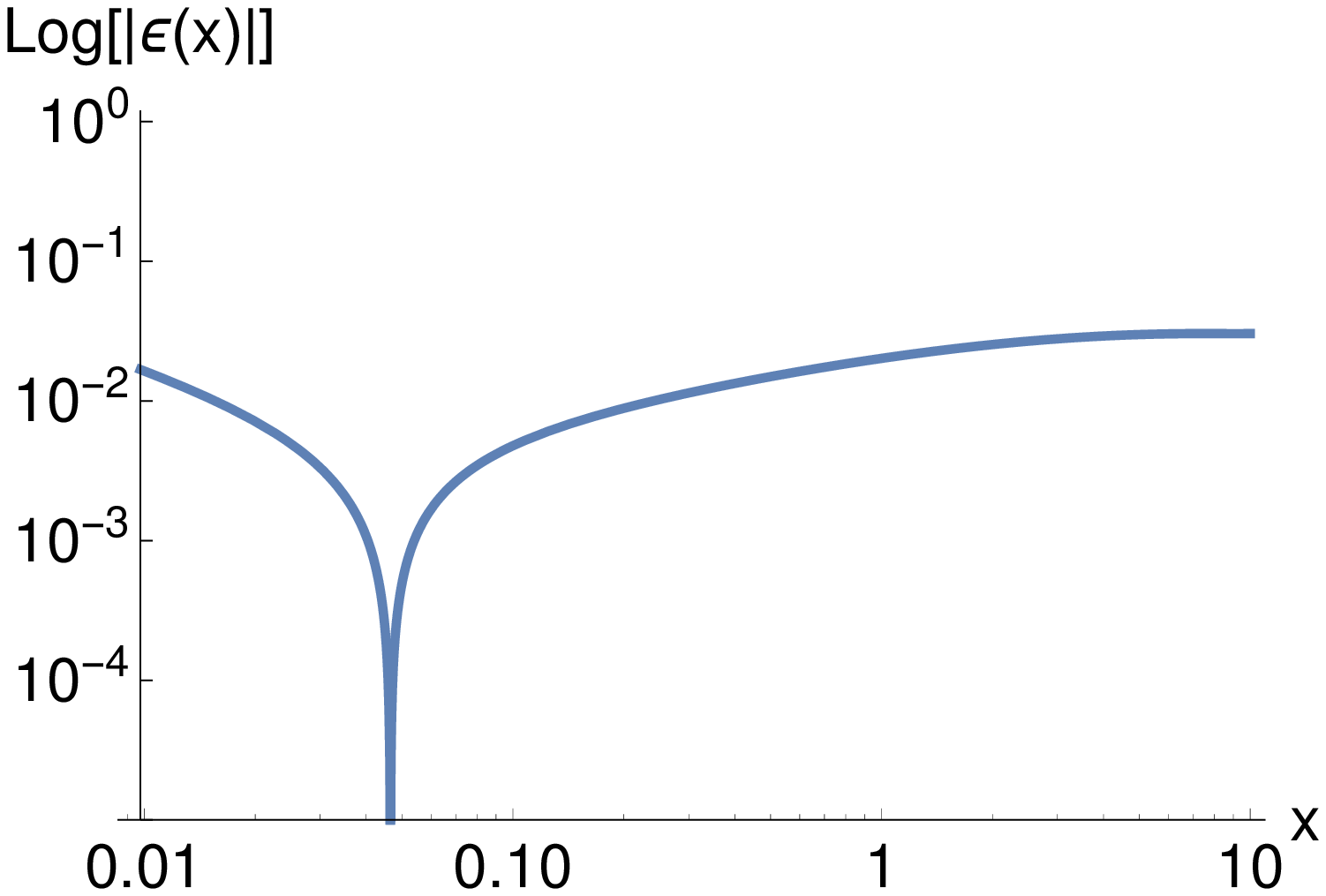}
 \caption{•}
\end{subfigure}%
\caption{
The application of the infinite integer derivative series for solving a linear fractional differential equation with constant coefficients. (a) Exact solution of the linear fractional differential equation, given in terms of the {\ML} function (solid red curve), is compared to the numerical solution (dotted blue curve) obtained by the fourth-order Runge-Kutta iterative method. The numerical solution is obtained by truncating the integer derivative expansion in Eq.(\ref{cap_result1}) for $q=1/2$ and retaining only the first three terms ($N=3$). (b) The log-log plot of the relative truncation error defined in Eq.(\ref{deferror1}) shows that truncating at $N=3$ results in sub one per-cent error.
}
\label{diffplot1}
\end{figure}

In the previous section we established convergence of the {\GL} fractional derivative by truncating the infinite integer derivative series and retaining only the first three terms. The goal of this section is to apply the truncated expansion of a fractional derivative to solve linear fractional differential equations (FDEs) with constant and variable coefficients. We choose two simplest non-trivial FDEs, which have solutions in terms of special functions, e.g.  Mittag-Leffler and generalized Fox-Wright functions. The comparison of the numerical approximation to the exact analytic result provides a direct test for the robustness of the numerical scheme based on the truncated expansion of a fractional derivative. 

The simplest form of the linear fractional differential equation with constant coefficients is given by 
\begin{equation}\label{diffequation1}
^{\mathrm{C}} \textbf{D}^q f(x) = -\lambda f(x),
\end{equation}
where $\lambda$ is a real-valued constant. The exact solution of Eq.(\ref{diffequation1}) is given in terms of the generalized {\ML} function \cite{kilbas2006theory} defined as
\begin{eqnarray}
E_{\alpha,\beta}(x) = \sum^{\infty}_{k=0}\frac{x^{k}}{\Gamma[\alpha k +\beta]}.
\end{eqnarray}
Specifically, the solution to Eq.(\ref{diffequation1}) is given by \cite{herrmann2014fractional}
\begin{eqnarray}
f(x) = E_{q}(-\lambda x^{q}) \equiv E_{q,1}(-\lambda x^q).
\end{eqnarray}
By adopting the Caputo fractional derivative, which ensures a solution convergent at the origin, and retaining the first $N=3$ terms in the integer derivative expansion, Eq.(\ref{cap_result1}),
for a fractional order $q=1/2$, 
we obtain a second order differential equation,
\begin{eqnarray}\label{difftruncated1}
-\frac{1}{6} x^2 f''(x)+x f'(x)+\sqrt{\pi x}  \lambda  f(x)+f(x)-f(0)=0.
\end{eqnarray}
The solution of the transformed differential equation is subject to the boundary conditions,
\begin{eqnarray}
f(0)=1\\\nn
\lim_{x\to\infty}f(x)=0.
\end{eqnarray}
The numerical solution of Eq.(\ref{difftruncated1}) is readily obtained via a 
fourth-order Runge-Kutta iterative method, shown in Fig.(\ref{diffplot1})
along with the relative truncation error $\epsilon(x)$ defined in Eq.(\ref{deferror1}).

Next we turn to a linear fractional differential equation with variable coefficients,
\begin{equation}\label{diffequation2}
^C\textbf{D}^{\alpha} f(x) = -\frac{\lambda f(x)}{x}.
\end{equation}
The exact solution to the fractional differential equation Eq.(\ref{diffequation2}) is given in terms of the generalized Fox-Wright function \cite{kilbas2006theory},
\begin{eqnarray}
{}_p\Psi_q \left[\begin{matrix} 
( a_1 , A_1 ) & ( a_2 , A_2 ) & \ldots & ( a_p , A_p ) \\ 
( b_1 , B_1 ) & ( b_2 , B_2 ) & \ldots & ( b_q , B_q ) \end{matrix} 
\;\bigg\rvert
z \right]
=
\sum_{n=0}^\infty 
\frac{\prod^{p}_{k=1}\Gamma( a_{k}+ A_{k} n )}
{\prod^{q}_{l=1}\Gamma( b_{l} + B_{l} n )} \, \frac {z^n} {n!}.
\end{eqnarray}
In particular, the solution to Eq.(\ref{diffequation2}) is given by 
\begin{eqnarray}
f(x)=
C x^{\alpha-1} 
{}_{0}\Psi_{1} \left[\begin{matrix} 
\--
 \\ 
(\alpha,\alpha-1) \end{matrix} 
\;\bigg\rvert
\frac{\lambda \hspace{0.5mm} x^{\alpha-1}}{1-\alpha}
 \right],
\end{eqnarray}
where $C$ is an arbitrary real constant.
In the special case of the fractional order $\alpha=1/2$ the generalized Fox-Wright function
is reduced to a Gaussian function,
\begin{eqnarray}
{}_{0}\Psi_{1} \left[\begin{matrix} 
\--
 \\ 
\left(\frac{1}{2},-\frac{1}{2}\right) \end{matrix} 
\;\bigg\rvert
z
 \right] = \frac{1}{\sqrt{\pi}}\exp(-z^{2}/4).
\end{eqnarray}
\begin{eqnarray}
{}_{0}\Psi_{1} \left[\begin{matrix} 
\--
 \\ 
\left(\frac{1}{2},-\frac{1}{2}\right) \end{matrix} 
\;\bigg\rvert
z
 \right] = \frac{1}{\sqrt{\pi}}\exp(-z^{2}/4)
\end{eqnarray}
where $``\--"$ in the argument of the Fox-Wright function stands for an absent argument.

\begin{figure}[t!]
\centering
\begin{subfigure}[b]{0.45\textwidth}
  \centering
  \captionsetup{justification=centering}
\includegraphics[width=\textwidth]{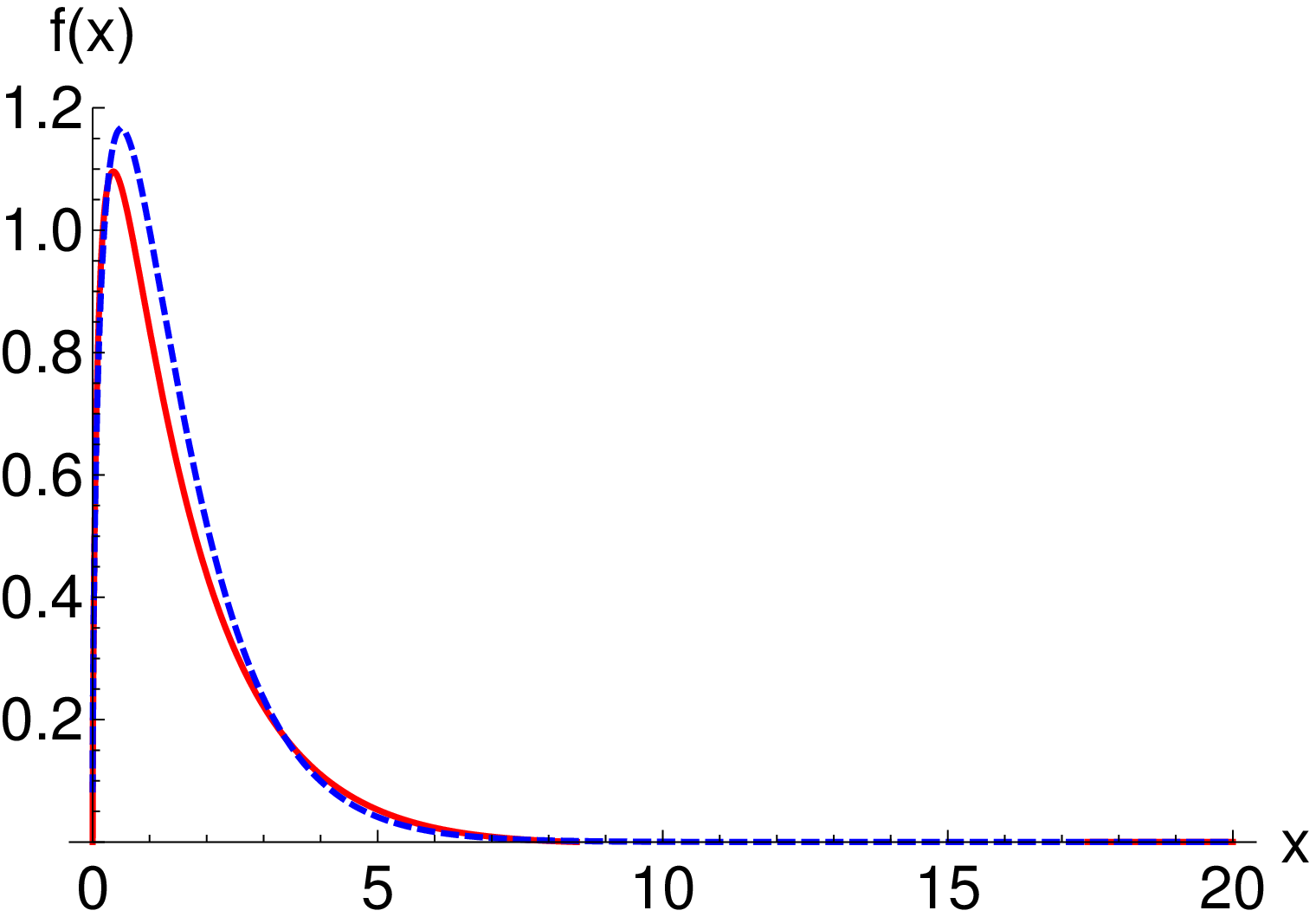}
 \caption{•}
\end{subfigure}%
\begin{subfigure}[b]{0.45\textwidth}
  \centering
  \captionsetup{justification=centering}
\includegraphics[width=\textwidth]{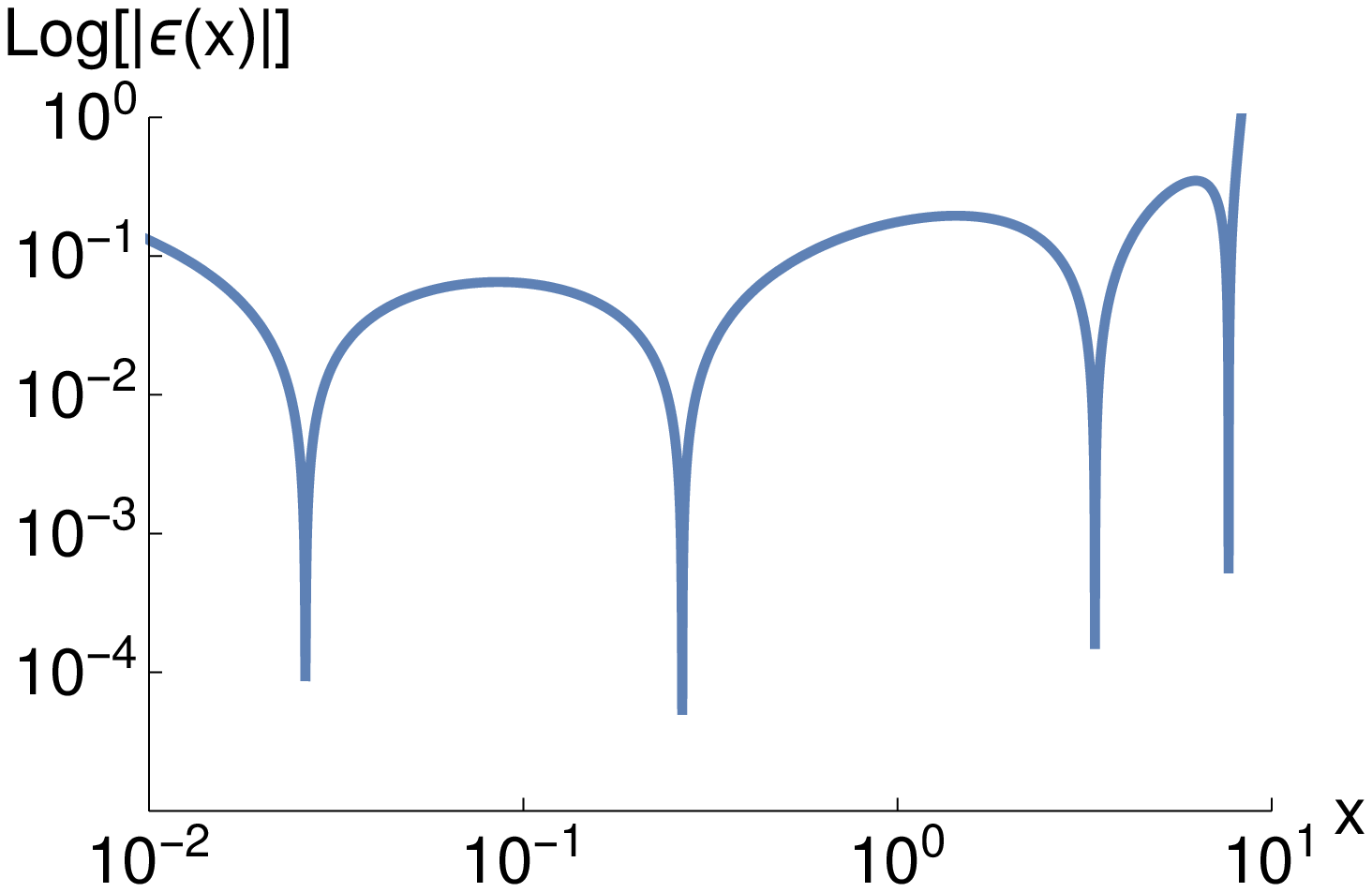}
 \caption{•}
\end{subfigure}%
\caption{
The application of the infinite integer derivative series for solving a linear fractional differential equation with variable coefficients. (a) The exact solution to the linear fractional differential equation with variable coefficients can be expressed in terms of the generalized Fox-Wright function  (solid red curve) which is compared to the numerical solution (dotted blue curve) obtained via a fourth-order Runge-Kutta iterative method. The numerical solution is obtained by truncating the integer derivative expansion in Eq.(\ref{cap_result1}) for $q=1/2$ and retaining only the first $N=3$ terms. (b) The log-log plot of the relative truncation error defined in Eq.(\ref{deferror1}) shows that truncating at $N=3$ results in sub ten per-cent error for $x \leq 10$.
}
\label{diffplot2}
\end{figure}

As a result the solution to the fractional differential equation Eq.(\ref{diffequation2}) in the special case of $\alpha=1/2$ is
\begin{equation}
f(x)= \frac{C}{\sqrt{\pi x}}\exp(-\lambda^{2}/x).
\end{equation}
If we further specify $f(1)=1$, we fix the constant $C$ and obtain the exact solution to Eq.(\ref{diffequation2}),
\begin{eqnarray}
f(x)= \frac{\exp(\lambda^{2}-\lambda^{2}/x)}{\sqrt{x}}.
\end{eqnarray}
By retaining the first $N=3$ terms in the integer derivative expansion Eq.(\ref{cap_result1}), we acquire,
\begin{eqnarray}
-\frac{1}{6} x^3 f''(x)+x^2 f'(x)+f(x) \left(x + \sqrt{\pi x } \lambda \right)-f(0) x=0.
\end{eqnarray}
If we change variables according to $x=1/y$, we obtain a transformed differential equation,
\begin{eqnarray}\label{diffequation3}
y^2 f''(y)+8 y f'(y)-6 (\sqrt{\pi y} + 1 ) f(y) = 0.
\end{eqnarray}
We specify the initial conditions as,
\begin{align}
f(0)&=0\\\nn
\left.\frac{\p f(y)}{\p x}\right|_{y=1}&=-\frac{1}{2},
\end{align}
and apply a fourth-order Runge-Kutta iterative method to find the numerical solution of 
Eq.(\ref{diffequation3}). The result along with the relative truncation error is shown in Fig.(\ref{diffplot2}). 

In this section, we successfully demonstrated that expanding a fractional derivative in terms of integer order derivatives is a robust method for solving linear fractional differential equations with both constant and variable coefficients. In the special case of a differential equation with variable coefficients, the truncated series with only the first $N=3$ terms leads to a $10\%$  error, while the very same truncation 
applied to a differential equation with constant coefficients results in a $1\%$ error. Although this method cannot be exhaustively tested for all possible fractional orders of differential operators and all types of FDEs, linear FDEs, considered in this work, constitute a large sample that can be used in many physical applications where the response of a system is proportional to a fractional order parameter \cite{hilfer2000applications, west2014colloquium, herrmann2014fractional}. Thus, the numerical scheme based on the truncated integer derivative expansion is a powerful method for solving a broad range of linear FDEs.
\cpg
\section{Conclusions}

In this paper we expressed the {\GL} fractional derivative as an infinite sum of integer order derivatives. We compared the obtained infinite expansion with the corresponding series produced by the {\RL} and {\Cp} definitions of a fractional derivative. We found that all three definitions are represented by the very same infinite series, with the exception of the lower index of summation for the {\Cp} fractional derivative which accounts for the initial conditions at the expansion point. Thus, we have shown that the integer derivative series representation provides a unified description for various definitions of a fractional derivative. 
  
By truncating the infinite expansion and retaining only the first few terms, we demonstrated the convergence of the {\GL} fractional derivative. We have shown that for functions represented by Taylor series with an infinite radius of convergence, the truncation error decreases with an increasing number of terms kept in the truncated expansion. We emphasized that the infinite expansion does not rely on the properties of the Taylor series, which has profound consequences for the functions characterized by a finite radius of convergence of the corresponding Taylor series. Specifically, we have shown that the infinite series of integer order derivatives for hyperbolic secant and tangent functions has an infinite radius of convergence, compared to the corresponding Taylor series with a finite radius of convergence of $\pi/2$. However, for a Gaussian function we found that the infinite expansion is divergent due to the factorial growth and oscillatory nature of the Hermite polynomials. Thus, the Gaussian function establishes limits of the universality of the infinite expansion of the {\GL} fractional derivative in terms of integer order derivatives.

Finally, we applied the truncated series for a fractional derivative to solve linear fractional differential equations with both constant and variable coefficients. We found that the fourth-order Runge-Kutta method applied to truncated fractional differential equations results in numerical solutions that rapidly converge to the exact solutions given in terms of Mittag-Leffler and generalized Fox-Wright special functions. 
Thus, we concluded that the integer derivative expansion can be adapted to a robust numerical method for solving linear fractional differential equations, such as the fractional {\Sch} and fractional diffusion equations.

\begin{acknowledgments}
The authors would like to thank David Benson, Daniel Jaschke, Nathan Smith, and Marc Valdez  for numerous and fruitful discussions. A.G., G.S. and L.D.C. acknowledge support from the US National Science Foundation under grant numbers PHY-1306638, PHY-1207881, PHY-1520915, and OAC-1740130, and the US Air Force Office of Scientific Research grant number FA9550-14-1-0287. This work was performed in part at the Aspen Center for Physics, which is supported by National Science Foundation grant PHY-1607611. U.A. acknowledges support from UAEU-UPAR(4) and UAEU-UPAR(7).
\end{acknowledgments}

\bibliography{math_fractional_lit1.bib}

%\nocite{*}
%\bibliography{aipsamp}% Produces the bibliography via BibTeX.

%\bibliography{ref}{}
%\bibliographystyle{plain}
%\begin{thebibliography}{99}
%\bibitem{c2} E. C. de Oliveira and J. A. T. Machado, ``A review of definitions for fractional derivatives and integrals,'' Mathematical Problems in Engineering \textbf{2014}, 238459 (2014). doi:10.1155/2014/238459
%\bibitem{c3} A. A. Kilbas, H. M. Srivastava, J. J. Trujillo, \textit{Theory and Applications of Fractional Differential Equations} (North-Holland Mathematics Studies, Vol. 204, Elsevier 2006) 
%\bibitem{c1} S. G. Samko, A. A. Kilbas, O. I. Marichev, \textit{Fractional Integrals and Derivatives} (Gorden and Breach Science Publishers 1993), p. 278
%\bibitem{c4} R. Herrmann, \textit{Fractional Calculus: An Introduction for Physicists} (World Scientific Publishing Co. Pte. Ltd. 2014), p. 411
%\bibitem{c5} J. J. Sakurai and J. Napolitano,\textit{ Modern Quantum Mechanics} (Addison-Wesley 2011), p. 46
%\end{thebibliography}

\end{document}